\documentclass[11pt,a4paper]{article}
\usepackage[utf8]{inputenc}
\usepackage{graphicx}
\newcounter{marnote}

\usepackage{setspace}
\usepackage{float}
\usepackage{subfig}
\usepackage{xr}
\usepackage{bm}
\usepackage[pdftex=true,hyperindex=true,colorlinks=true]{hyperref}
\usepackage{mathtools}
\usepackage{indentfirst}
\usepackage[T1]{fontenc}
\usepackage{color}
\usepackage{amsfonts} 
\usepackage{amsmath}

\usepackage{mathrsfs}
\usepackage{tikz}
\usepackage{dsfont}
\usetikzlibrary{calc}
\usetikzlibrary{patterns}
\usepackage{pstricks}
\usepackage{esint}
\tikzset{math3d/.style={x={(-0.353cm,-0.353cm)},z={(0cm,1cm)},y={(1cm,0cm)}}}

\newtheorem{Theo}{Theorem}[section] 
\newtheorem{Prop}{Proposition}[section]
\newtheorem{Not}{Notation}[section]
\newtheorem{Lem}{Lemma}[section]

\newtheorem{Def}{Definition}[section]
\newtheorem{Rem}{Remark}[section]

\newenvironment{Dem}[1]
{\begin{trivlist}\item[\hskip5mm{\sc #1}]}{\end{trivlist}}

\DeclareMathOperator{\divg}{div}
\newcommand{\R}{\mathbb{R}}
\newcommand{\Z}{\mathbb{Z}}

\newcommand{\dx}{\ \textrm{d}x}

\newcommand{\harpoon}[1]{-\hspace{-0.2cm}\rightharpoonup}
\newcommand{\qed}{\hfill \ensuremath{\Box}}

\def\Sum{\displaystyle\sum}

\newpage
\usepackage{tikz}
\usetikzlibrary{patterns}
\usetikzlibrary{arrows}

\tikzset { xmin/.store in=\xmin, xmin/.default=-3, xmin=-3,
xmax/.store in=\xmax, xmax/.default=3, xmax=3,
ymin/.store in=\ymin, ymin/.default=-3, ymin=-3,
ymax/.store in=\ymax, ymax/.default=3, ymax=3}

\tikzset {domaine/.style 2 args={domain=#1:#2}}

\makeatletter
\renewcommand{\thefigure}{\ifnum \c@section>\z@ \thesection.\fi
\@arabic\c@figure}
\@addtoreset{figure}{section}
\makeatother
\makeatletter
\renewcommand{\theequation}{\ifnum \c@section>\z@ \thesection.\fi
\@arabic\c@equation}
\@addtoreset{equation}{section}
\makeatother

\title{\LARGE Homogenization of high-contrast two-phase conductivities perturbed by a magnetic field. Comparison between dimension two and dimension three.}
\author{\begin{tabular}{cc}
    {Marc BRIANE } & {Laurent PATER}
    \\*[-.0em]
    {\small  Institut de Recherche Math\'ematique de Rennes} & {\small Institut de Recherche Math\'ematique de Rennes}
    \\*[-.3em]
    {\small  INSA de Rennes} &	 {\small Universit\'e de Rennes 1}
    \\*[-.3em]
    {\small mbriane@insa-rennes.fr} & {\small laurent.pater@ens-cachan.org}
\end{tabular}}
\begin{document}
\maketitle
\textcolor{white}{mon texte}
\begin{abstract}
Homogenized laws for sequences of high-contrast two-phase non-symmetric conductivities perturbed by a parameter $h$ are derived in two and three dimensions. The parameter $h$ characterizes the antisymmetric part of the conductivity for an idealized model of a conductor in the presence of a magnetic field. In dimension two an extension of the Dykhne transformation to non-periodic high conductivities permits to prove that the homogenized conductivity depends on $h$ through some homogenized matrix-valued function obtained in the absence of a magnetic field. This result is improved in the periodic framework thanks to an alternative approach, and illustrated by a cross-like thin structure. Using other tools, a fiber-reinforced medium in dimension three provides a quite different homogenized conductivity.
\end{abstract}
\section{Introduction}
The mathematical theory of homogenization for second-order elliptic partial differential equations has been widely studied since the pioneer works of Spagnolo on $G$-convergence \cite{Spa}, of Murat, Tartar on $H$-convergence \cite{Mur1,Mur2}, and of Bensoussan, Lions, Papanicolaou on periodic structures \cite{BLP}, in the framework of uniformly bounded (both from below and above) sequences of conductivity matrix-valued functions. It is also known since the end of the seventies \cite{FeKh,Khr} (see also the extensions \cite{BelBou,BrTc,Bri1,KhMa}) that the homogenization of the sequence of conductivity problems, in a bounded open set $\Omega$ of $\R^3$,
\begin{equation}
\left\{\begin{array}{rl}
\mathrm{div}\left(\sigma_n \nabla u_n\right)= f & \mbox{in }\Omega
\\
u_n = 0 & \mbox{on }\partial\Omega,
\end{array}\right.
\label{pbintro}
\end{equation}
with a uniform boundedness from below but not from above for $\sigma_n$, may induce nonlocal effects. However, the situation is radically different in dimension two since the nature of problem (\ref{pbintro}) is shown \cite{Brinonloc, bcconvunif} to be preserved in the homogenization process provided that the sequence $\sigma_n$ is uniformly bounded from below.
\par
H-convergence theory includes the case of non-symmetric conductivities in connection with the Hall effect \cite{Hall} in electrodynamics (see, e.g., \cite{LaLi,Ali}). Indeed, in the presence of a constant magnetic field the conductivity matrix is modified and becomes non-symmetric. Here, we consider an idealized model of an isotropic conductivity $\sigma(h)$ depending on a parameter $h$ which characterizes the antisymmetric part of the conductivity in the following way:\begin{itemize} \item in dimension two, \begin{equation}
\sigma(h) = \alpha I_2 + \beta h J,
\quad J:=\left(\begin{smallmatrix} 0 & -1 \\ 1 & 0 \end{smallmatrix}\right),\label{model2}
\end{equation}
where $\alpha$, $\beta$ are scalar an $h \in \mathbb R$,
\item in dimension three,\begin{equation}
\sigma(h) = \alpha I_3 + \beta \mathscr{E}(h),
\quad\mathscr{E}(h):=\left(\begin{smallmatrix}
0    & -h_3  & h_2 \cr
h_3 & 0    & -h_1 \cr
-h_2 & h_1 & 0 \cr
\end{smallmatrix}\right),\label{model3}
\end{equation}
where $\alpha$, $\beta$ are scalar and $h \in \mathbb R^3$.
\end{itemize} Since the seminal work of Bergman \cite{Ber} the influence of a low magnetic field in composites has been studied for two-dimensional composites \cite{MilHall,BerLiSt,BMM}, and for columnar composites \cite{BerSt1,BerSt2,StrelBer,Grab1,Grab2}. The case of a strong field, namely when the symmetric part and the antisymmetric part of the conductivity are of the same order, has been also investigated \cite{BerSt1a,BerSS}. Moreover, dimension three may induce anomalous homogenized Hall effects \cite{BM1,BM2,BM3} which do not appear in dimension two~\cite{BMM}.
\par
In the context of high-contrast problems the situation is more delicate when the conductivities are not symmetric. An extension in dimension two of H-convergence for non-symmetric and non-uniformly bounded conductivities was obtained in \cite{BrCa1} thanks to an appropriate div-curl lemma. More recently, the Keller, Dykhne \cite{Keller,Dyk} two-dimensional duality principle which claims that the mapping
\begin{equation}
A\mapsto{\frac{A^\text{T}}{ \det A}}
\label{dual}
\end{equation}
is stable under homogenization, was extended to high-contrast conductivities in \cite{NHM}. However, the homogenization of both high-contrast and non-symmetric conductivities has not been precisely studied in the context of the strong field magneto-transport especially in dimension three. In this paper we establish an effective perturbation law for a mixture of two high-contrast isotropic phases in the presence of a magnetic field. The two-dimensional case is performed in a general way for non-periodic and periodic microstructures. It is then compared to the case of a three-dimensional fiber-reinforced microstructure.

\bigskip
In dimension two, following the modelization \eqref{model2}, consider a sequence $\sigma_n(h)$ of isotropic two-phase matrix-valued conductivities perturbed by a fixed constant $h\in\R$, and defined by
\begin{equation}
\sigma_n(h):=(1-\chi_n)\big(\alpha_1 I_2 + \beta_1h\,J\big)+\chi_n\big(\alpha_{2,n} I_2+\beta_{2,n}h\,J\big),
\label{intro0}
\end{equation}
where $\chi_n$ is the characteristic function of phase $2$, with volume fraction $\theta_n\to 0$, $\alpha_1 >0$, $\beta_1$ are the constants of the low conducting phase~$1$, and $\alpha_{2,n} \to \infty$, $\beta_{2,n}$ are real sequences of the highly conducting phase~$2$ where $\beta_{2,n}$ is possibly unbounded. The coefficients $\alpha_1$ and $\beta_1$, respectively $\alpha_{2,n}$ and $\beta_{2,n}$ also have the same order of magnitude according to the strong field assumption. Assuming that the sequence $\theta_n^{-1}\chi_n$ converges weakly-$*$ in the sense of the Radon measures to a bounded function, and that $\theta_n\alpha_{2,n}$, $\theta_n\beta_{2,n}$ converge respectively to constants $\alpha_2 >0$, $\beta_2$, we prove (see Theorem~\ref{divcurl1}) that the perturbed conductivity $\sigma_n(h)$ converges in an appropriate sense of H-convergence (see Definition~\ref{H(M)}) to the homogenized matrix-valued function
\begin{equation}
\sigma_*(h)=\sigma_*^0\big(\alpha_1, \alpha_2 + \alpha_2^{-1}\beta_2^2\,h^2\big)+\beta_1h\,J,
\label{intro1}
\end{equation}
for some matrix-valued function $\sigma_*^0$ which depends uniquely on the microstructure $\chi_n$ in the absence of a magnetic field, and is defined for a subsequence of $n$. 
The proof of the result is based on a Dykhne transformation of the type
\begin{equation}
A_n\mapsto\big((p_n A_n + q_n J)^{-1} + r_n J\big)^{-1},
\label{Dtrans}
\end{equation}
which permits to change the non-symmetric conductivity $\sigma_n(h)$ into a symmetric one. Then, extending the duality principle (\ref{dual}) established in \cite{NHM}, we prove that transformation (\ref{Dtrans}) is also stable under high-contrast conductivity homogenization.
\par
In the periodic case, i.e. when $\sigma_n(h)(\cdot)=\Sigma_n(\cdot/\varepsilon_n)$ with $\Sigma_n$ $Y$-periodic and $\varepsilon_n\to 0$, we use an alternative approach based on an extension of Theorem~4.1 of \cite{bcconvunif} to $\varepsilon_n Y$-periodic but non-symmetric conductivities (see Theorem \ref{thunifconv}). So, it turns out that the homogenized conductivity $\sigma_*(h)$ is the limit as $n\to\infty$ of the constant H-limit $(\sigma_n)_*$ associated with the periodic homogenization (see, e.g., \cite{BLP}) of the oscillating sequence $\Sigma_n(\cdot/\varepsilon)$ as $\varepsilon\to 0$ and for a fixed $n$. Finally, the Dykhne transformation performed by Milton \cite{MilHall} (see also \cite{MilBook}, Chapter~4) applied to the local periodic conductivity $\Sigma_n$ and its effective conductivity $(\sigma_n)_*$, allows us to recover the perturbed homogenized formula (\ref{intro1}).
An example of a periodic cross-like thin structure provides an explicit computation of $\sigma_*(h)$ (see Proposition~\ref{thcross}).
\par
To make a comparison with dimension three we restrict ourselves to the $\varepsilon_n Y$-periodic fiber-reinforced structure introduced by Fenchenko, Khruslov \cite{FeKh} to derive a nonlocal effect in homogenization. However, in the present context the fiber radius $r_n$ is chosen to be super-critical, i.e. $r_n\to 0$ and $\varepsilon_n^2 |\ln r_n|\to 0$, in order to avoid such an effect. Similarly to (\ref{intro0}) and following the modelization \eqref{model3}, the perturbed conductivity is defined for $h\in\R^3$, by
\begin{equation}
\sigma_n(h):=(1-\chi_n)\big(\alpha_1 I_3 + \beta_1\,\mathscr{E}(h)\big)+\chi_n\big(\alpha_{2,n} I_3+\beta_{2,n}\,\mathscr{E}(h)\big),
\label{intro3}
\end{equation}
where $\chi_n$ is the characteristic function of the fibers which are parallel to the direction $e_3$.
The form of \eqref{intro3} ensures the rotational invariance of $\sigma_n(h)$ for those orthogonal transformations which leave $h$ invariant. Under the same assumptions on the conductivity coefficients as in the two-dimensional case, with $\theta_n=\pi\,r_n^2$, but using a quite different approach, the homogenized conductivity is given by (see Theorem \ref{thdim3})
\begin{equation}\label{introo}
\sigma_*(h)=\alpha_1 I_3 + \displaystyle \left(\frac{\alpha_2^3 + \alpha_2 \beta_2^2 |h|^2}{\alpha_2^2 + \beta_2^2 h_3^2} \right) e_3 \otimes e_3 + \beta_1 \mathscr{E}(h).
\end{equation}
The difference between formulas \eqref{intro1} and \eqref{introo} provides a new example of gap between dimension two and dimension three in the high-contrast homogenization framework. As former examples of dimensional gap, we refer to the works \cite{BMM,BM1} about the 2d positivity property, versus the 3d non-positivity, of the effective Hall coefficient, and to the works \cite{bcconvunif,FeKh} concerning the 2d lack, versus the 3d appearance, of nonlocal effects in the homogenization process.

\bigskip
The paper is organized as follows. Section~\ref{nonper} and~\ref{per} deal with dimension two. In Section~\ref{nonper} we study the two-dimensional general (non-periodic) case thanks to an appropriate div-curl lemma. In Section~\ref{per} an alternative approach is performed in the periodic framework. Finally, Section~\ref{dim3} is devoted to the three-dimensional case with the fiber-reinforced structure.

\subsection*{Notations}

\begin{itemize}
\item[$\bullet$] $\Omega$ denotes a bounded open subset of $\R^d$;
\item[$\bullet$] $I_d$ denotes the unit matrix in $\R^{d \times d}$ , and $J:=\left(\begin{smallmatrix} 0 & -1 \\ 1 & 0 \end{smallmatrix}\right)$;
\item[$\bullet$] for any matrix $A$ in $\R^{d \times d}$, $A^\mathrm{T}$ denotes the transposed of the matrix $A$, $A^s$ denotes its symmetric part;
\item[$\bullet$] for $h \in \R^3$, $\mathscr E(h)$ denotes the antisymmetric matrix in $\R^{3 \times 3}$ defined by $\mathscr E(h)\,x:=h \times x$, for $x\in\R^3$;
\item[$\bullet$] for any $A,B\in\R^{d \times d}$, $A \leq B$ means that for any $\xi \in \R^d$, $A \xi \cdot \xi \leq B \xi \cdot \xi$; we will use the fact that for any invertible matrix $A\in\R^{d \times d}$, $A \geq \alpha I_d\Rightarrow A^{-1} \leq \alpha^{-1} I_d$;
\item[$\bullet$] $| \cdot |$ denotes both the euclidean norm in $\R^d$ and the subordinate norm in $\R^{d \times d}$;
\item[$\bullet$] for any locally compact subset $X$ of $\R^d$, $\mathcal M(X)$ denotes the space of the Radon
measures defined on $X$;
\item[$\bullet$] for any $\alpha, \beta >0$, $\mathcal M(\alpha, \beta; \Omega)$ is the set of the invertible matrix-valued functions $A : \Omega\to\R^{d \times d}$ such that \begin{equation}
\forall\,\xi \in \R^d, \quad A(x) \xi \cdot \xi \geq \alpha |\xi|^2 \quad \mbox{and} \quad A^{-1}(x) \xi \cdot \xi \geq \beta^{-1} |\xi|^2 \quad \mbox{a.e. in } \Omega;
\label{M(a,b)}
\end{equation}
\item[$\bullet$] $C$ denotes a constant which may vary from a line to another one.
\end{itemize}

In the sequel, we will use the following extension of $H$-convergence and introduced in \cite{NHM}:

\begin{Def} Let $\alpha_n$ and $\beta_n$ be two sequences of positive numbers such that $\alpha_n \leq \beta_n$, and let $A_n$ be a sequence of matrix-valued functions in $\mathcal M(\alpha_n,\beta_n;\Omega)$
(see \eqref{M(a,b)}).

The sequence $A_n$ is said to $H(\mathcal M(\Omega)^2 )$-converge to the matrix-valued function $A_*$ if for any distribution $f$ in $H^{-1} (\Omega)$, the solution $u_n$ of the problem
\begin{equation*}
\left\{\begin{array}{rl}
\mathrm{div}\left(A_n \nabla u_n\right)= f & \mbox{in }\Omega
\\
u_n = 0 & \mbox{on }\partial\Omega,
\end{array}\right.
\end{equation*}satisfies the convergences
\begin{equation*}
\left\{\begin{array}{rll}
u_n & \harpoon 0 \ u & \mbox{in }H^1_0(\Omega)
\\
A_n \nabla u_n & \harpoon 0 \ A_* \nabla u & \mbox{weakly-$*$ in }\mathcal M(\Omega)^2,
\end{array}\right.
\end{equation*}where $u$ is the solution of the problem\begin{equation*}
\left\{\begin{array}{rl}
\mathrm{div}\left(A_* \nabla u\right)= f & \mbox{in }\Omega
\\
u = 0 & \mbox{on }\partial\Omega.
\end{array}\right.
\end{equation*}
\label{H(M)}
\end{Def}

We now give a notation for $H(\mathcal M(\Omega)^2 )$-limits of high-contrast two-phase composites. We consider the characteristic function $\chi_n$ of the highly conducting phase, and denote $\omega_n:= \{ \chi_n = 1 \}$.

\begin{Not} A sequence of isotropic two-phase conductivities in the absence of a magnetic field is denoted by
\begin{equation}
\sigma_n^0(\alpha_{1,n},\alpha_{2,n}):= (1-\chi_n) \alpha_{1,n} I_2 + \chi_n \alpha_{2,n} I_2,
\end{equation}with\begin{equation}
\lim \limits_{n \to \infty} \alpha_{1,n} = \alpha_1 > 0 \quad \text{and} \quad  \lim \limits_{n \to \infty} \left| \omega_n \right| \alpha_{2,n} = \alpha_2 > 0, \end{equation} and its $H(\mathcal M(\Omega)^2)$-limit is denoted by $\sigma_*^0(\alpha_1,\alpha_2)$.
\label{defabs}
\end{Not}

\section{A two-dimensional non-periodic medium\label{nonper}}

\subsection{A div-curl approach}

We extend the classical div-curl lemma.

\begin{Lem}

Let $\Omega$ be a bounded open subset of $\R^2$. Let $\alpha > 0$, let $\bar{a} \in L^\infty(\Omega)$ and let $A_n$ be a sequence of matrix-valued functions in $L^\infty(\Omega)^{2 \times 2}$ (not necessarily symmetric) satisfying\begin{equation}A_n \geq \alpha I_2 \ \ \textrm{a.e.} \ \textrm{in} \ \Omega \ \ \ \ \textrm{and} \ \ \ \ \displaystyle{\frac{\det A_n}{\det A_n^s}} |A_n^s| \rightharpoonup \bar{a} \in L^\infty(\Omega) \ \ \textrm{weakly-$*$} \ \textrm{in} \ \mathcal{M}(\Omega).\label{conddet1}\end{equation}

\noindent Let $\xi_n$ be a sequence in $L^2(\Omega)^2$ and $v_n$ a sequence in $H^1(\Omega)$ satisfying the following assumptions:

(i) $\xi_n$ and $v_n$ satisfy the estimate\begin{equation}\int_\Omega A_n^{-1} \xi_n \cdot \xi_n \textrm{d}x + ||v_n||_{H^1(\Omega)} \leq C ; \label{divrot1}\end{equation}

(ii) $\xi_n$ satisfies the classical condition\begin{equation}\divg \xi_n \ \textrm{is} \ \textrm{compact} \ \textrm{in} \ H^{-1}(\Omega)\label{divrot2}.\end{equation}

\noindent Then, there exist $\xi$ in $L^2(\Omega)^2$ and $v$ in $H^1(\Omega)$ such that the following convergences hold true up to a subsequence\begin{equation}\xi_n \rightharpoonup \xi \ \ \textrm{weakly-$*$} \ \textrm{in} \ \mathcal{M}(\Omega)^2 \ \ \ \ \textrm{and} \ \ \ \ \nabla v_n \rightharpoonup \nabla v \quad \textrm{weakly} \ \textrm{in} \ L^2(\Omega)^2.\label{convdivcurl}\end{equation}

\noindent Moreover, we have the following convergence in the distribution sense$$\xi_n \cdot \nabla v_n \rightharpoonup \xi \cdot \nabla v \ \ \textrm{weakly} \ \textrm{in} \ \mathscr{D}'(\Omega).$$
\label{divcurl}
\end{Lem}
\noindent \textbf{Proof of Lemma \ref{divcurl}}.
The proof consists in considering the "good-divergence" sequence $\xi_n$ as a sum of a compact sequence of gradients $\nabla u_n$ and a sequence of divergence-free functions $J\nabla z_n$. We then use Lemma 3.1 of \cite{NHM} to obtain the strong convergence of $z_n$ in $L^2_{loc}(\Omega)$. Finally, replacing $\xi_n$ by $\nabla u_n + J \nabla z_n$, we conclude owing to integration by parts.

\bigskip

\noindent \textit{First step}: Proof of convergences \eqref{convdivcurl}.

\smallskip

An easy computation gives \begin{equation}
\Big(\big(A_n^{-1}\big)^s\Big)^{-1} = \displaystyle{\frac{\det A_n}{\det A_n^s}} A_n^s. \label{calceasy}
\end{equation}

The sequence $\xi_n$ is bounded in $L^1(\Omega)^{2}$ since the Cauchy-Schwarz inequality combined with the weak-$*$ convergence of \eqref{conddet1}, \eqref{divrot1} and \eqref{calceasy} yields \begin{equation*}
\left(\int_\Omega |\xi_n| \dx \right)^2 \leq \int_{\Omega} \left| \Big(\big(A_n^{-1}\big)^s\Big)^{-1} \right| \dx \int_{\Omega} \big( A_n^{-1} \big)^s \xi_n \cdot \xi_n \dx = \int_{\Omega} \displaystyle{\frac{\det A_n}{\det A_n^s}} \left| A_n^s \right| \dx \int_{\Omega} A_n^{-1} \xi_n \cdot \xi_n \dx \leq C.
\end{equation*}

\noindent Therefore, $\xi_n$ converges up to a subsequence to some $\xi \in \mathcal{M}(\Omega)^2$ in the weak-$*$ sense of the measures. Let us prove that the vector-valued measure $\xi$ is actually in $L^2 (\Omega)^2$. Again by the Cauchy-Schwarz inequality combined with \eqref{conddet1}, \eqref{divrot1} and \eqref{calceasy} we have, for any $\Phi \in \mathscr C_0(\Omega)^2$,

\begin{align*}
&\left|\int_{\Omega} \xi(\mathrm{d}x) \cdot \Phi \right| = \lim \limits_{n \to \infty} \left|\int_{\Omega} \xi_n \cdot \Phi \dx \right| \\
& \leq \limsup \limits_{n \to \infty} \left(\int_{\Omega} \displaystyle{\frac{\det A_n}{\det A_n^s}} \left| A_n^s \right| |\Phi|^2 \dx \right)^{\frac{1}{2}} \left(\int_{\Omega} A_n^{-1} \xi_n \cdot \xi_n \dx \right)^{\frac{1}{2}} \leq C \left(\int_{\Omega} \bar{a} |\Phi|^2 \dx \right)^{\frac{1}{2}},
\end{align*} which implies that $\xi$ is absolutely continuous with respect to the Lebesgue measure. Since $\bar{a} \in L^\infty(\Omega)$, we also get that
\begin{equation*}
\left|\int_{\Omega} \xi \cdot \Phi \dx \right| \leq ||\Phi||_{L^2(\Omega)^2}
\end{equation*}hence $\xi \in L^2 (\Omega)^2$. Therefore, the first convergence of \eqref{convdivcurl} holds true with its limit in $L^2 (\Omega)^2$. The second one is immediate.

\bigskip

\noindent \textit{Second step}: Introduction of a stream function.

\smallskip

By \eqref{divrot2}, the sequence $u_n$ in $H^1_0(\Omega)$ defined by $u_n:= \Delta^{-1} \left( \mathrm{div} \ \xi_n \right)$ strongly converges in $H^1_0(\Omega)$:\begin{equation}
u_n \longrightarrow u \quad \text{in} \ H^1_0(\Omega). \label{unconv}
\end{equation}
Let $\omega$ be a regular simply connected open set such that $\omega \subset \subset \Omega$. Since by definition $\xi_n - \nabla u_n$ is a divergence-free function in $L^2(\Omega)^2$, there exists (see, e.g., \cite{GirRav}) a unique stream function $z_n \in H^1(\omega)$ with zero $\omega$-average such that \begin{equation} \label{egalstream}
\xi_n = \nabla u_n + J \nabla z_n \quad \text{a.e. in} \ \omega.
\end{equation}

\bigskip

\noindent \textit{Third step}: Convergence of the stream function $z_n$.

\smallskip

Since $\nabla u_n$ is bounded in $L^2(\Omega)^2$ by the second step, $\xi_n$ is bounded in $L^1 (\Omega)^2$ by the first step and $z_n$ has a zero $\omega$-average, the Sobolev embedding of $W^{1,1}(\omega)$ into $L^2(\omega)$ combined with the Poincaré-Wirtinger inequality in $\omega$ implies that $z_n$ is bounded in $L^2(\omega)$ and thus converges, up to a subsequence still denoted by $n$, to a function $z$ in $L^2(\omega)$.

\noindent Moreover, let us define \[S_n:= \left( J^{-1} (A_n^{-1})^s J \right)^{-1}.\]The Cauchy-Schwarz inequality gives\begin{align*}
&\int_\omega S_n^{-1} \nabla z_n \cdot\nabla z_n \dx = \int_\omega J^{-1} (A_n^{-1})^s J \nabla z_n \cdot\nabla z_n \dx \\
&= \int_\omega (A_n^{-1})^s J \nabla z_n \cdot J \nabla z_n \dx \\
& = \int_\omega (A_n^{-1})^s \big[\xi_n - \nabla u_n\big] \cdot \big[\xi_n - \nabla u_n] \dx \\
& \leq 2 \int_\omega (A_n^{-1})^s \xi_n \cdot \xi_n \dx + 2 \int_\omega (A_n^{-1})^s \nabla u_n \cdot \nabla u_n \dx \\
& = 2 \int_\omega A_n^{-1} \xi_n \cdot \xi_n \dx + 2 \int_\omega A_n^{-1} \nabla u_n \cdot \nabla u_n \dx .\end{align*}

\noindent The first term is bounded by \eqref{divrot1} and the last term by the inequality
\noindent $A_n^{-1} \leq \alpha^{-1} I_2$ and the convergence \eqref{unconv}. Therefore, the sequences $v_n:= z_n$ and, by \eqref{conddet}, $S_n$ satisfy all the assumptions of Lemma 3.1 of \cite{NHM} since, by \eqref{calceasy}, \[ S_n = \displaystyle{\frac{\det A_n}{\det A_n^s}} J^{-1} A_n^s J .\] Then, we obtain the convergence \begin{equation}z_n \longrightarrow z \quad \mathrm {strongly} \ \mathrm{in} \ L_{\textrm{loc}}^2(\omega). \label{z}\end{equation} Moreover, the convergence \eqref{unconv} gives \begin{equation}
\xi = \nabla u + J \nabla z \quad \text{in} \ \mathscr{D}'(\omega).\label{xiu}
\end{equation}

\bigskip

\noindent \textit{Fourth step}: Integration by parts and conclusion.

\smallskip

We have, as $J \nabla v_n$ is a divergence-free function, \begin{equation} \label{conc}
\xi_n \cdot \nabla v_n = \left( \nabla u_n + J \nabla z_n \right) \cdot \nabla v_n = \nabla u_n \cdot \nabla v_n - \mathrm{div} \left( z_n J \nabla v_n \right).
\end{equation}The strong convergence of $\nabla u_n$ in \eqref{unconv}, the second weak convergence of \eqref{convdivcurl} justified in the first step and \eqref{z} give \begin{equation}
\label{finaldivcurl} \nabla u_n \cdot \nabla v_n - \mathrm{div} \left( z_n J \nabla v_n \right) \ \harpoon 0 \nabla u \cdot \nabla v - \mathrm{div} \left( z J \nabla v \right)\quad \mathrm{in} \ \mathscr D'(\omega).
\end{equation}We conclude, by combining this convergence with \eqref{conc}, \eqref{xiu} and integrating by parts, to the convergence $$\xi_n \cdot \nabla v_n \ \harpoon 0 \nabla u \cdot \nabla v - \mathrm{div} \left( z J \nabla v \right) =  (\nabla u + J \nabla z) \cdot \nabla v = \xi \cdot \nabla v \ \ \textrm{weakly} \ \textrm{in} \ \mathscr{D}'(\omega).$$for an arbitrary open subset $\omega$ of $\Omega$.\qed

\bigskip

For the reader's convenience, we first recall in Theorem \ref{theo22} below the main result of \cite{NHM} concerning the Keller duality for high contrast conductivities. Then, Proposition \ref{lemDykhnelim} is an extension of this result to a more general transformation.

\begin{Theo}[\cite{NHM}]

Let $\Omega$ be a bounded open subset of $\R^2$ such that $|\partial \Omega| = 0 $. Let $\alpha > 0$, let $\beta_n$, $n \in \mathbb N$ be a sequence of real numbers such that $\beta_n \geq \alpha$, and let $A_n$ be a sequence of matrix-valued functions (not necessarily symmetric) in $\mathcal M(\alpha,\beta_n ; \Omega)$. Assume that there exists a function $\bar{a} \in L^\infty(\Omega)$ such that\begin{equation}\displaystyle{\frac{\det A_n}{\det A_n^s}} |A_n^s| \rightharpoonup \bar{a} \ \ \textrm{weakly-$*$} \ \textrm{in} \ \mathcal{M}(\Omega).\end{equation}

\noindent Then, there exist a subsequence of $n$, still denoted by $n$, and a matrix-valued function $A_*$ in \linebreak[4] $\mathcal M(\alpha, \beta ; \Omega)$, with $\beta = 2 ||\bar{a}||_{L^\infty(\Omega)}$, such that\begin{equation}A_n \stackrel{H(\mathcal{M}(\Omega)^2)}{-\hspace{-0.2cm}\rightharpoonup}{A_*} \ \ \ \ \textrm{and} \ \ \ \ \frac{A_n^\text{T}}{\det A_n} \stackrel{H(\mathcal{M}(\Omega)^2)}{-\hspace{-0.2cm}\rightharpoonup}{\frac{A_*^\text{T}}{\det A_*}}.\end{equation}\label{theo22}

\end{Theo}

\begin{Prop}

Let $\Omega$ be a bounded open subset of $\R^2$ such that $|\partial \Omega| = 0 $. Let $p_n$, $q_n$ and $r_n$, $n \in \mathbb N$ be sequences of real numbers converging respectively to $p > 0$, $q$ and $0$. Let $\alpha > 0$, let $\beta_n$, $n \in \mathbb N$ be a sequence of real numbers such that $\beta_n \geq \alpha$, and let $A_n$ be a sequence of matrix-valued functions in $\mathcal M(\alpha,\beta_n ; \Omega)$ (not necessarily symmetric) satisfying\begin{equation}r_n A_n \ \textrm{is} \ \textrm{bounded} \ \textrm{in} \ L^\infty(\Omega)^{2 \times 2} \ \ \ \ \textrm{and} \ \ \ \ \displaystyle{\frac{\det A_n}{\det A_n^s}} |A_n^s| \rightharpoonup \bar{a} \in L^\infty(\Omega) \ \ \textrm{weakly-$*$} \ \textrm{in} \ \mathcal{M}(\Omega), \label{conddet}\end{equation}

\noindent and that \begin{equation} B_n =\big((p_n A_n + q_n J)^{-1} + r_n J\big)^{-1} \ \ \textrm{is} \ \textrm{a} \ \textrm{sequence} \ \textrm{of} \ \textrm{symmetric} \ \textrm{matrices}.\end{equation}

\noindent Then, there exist a subsequence of $n$, still denoted by $n$, and a matrix-valued function $A_*$ in \linebreak[4] $\mathcal M(\alpha, \beta ; \Omega)$, with $\beta = 2 ||\bar{a}||_{L^\infty(\Omega)}$, such that\begin{equation}A_n \stackrel{H(\mathcal{M}(\Omega)^2)}{-\hspace{-0.2cm}\rightharpoonup}{A_*} \ \ \ \ \textrm{and} \ \ \ \ \big( ( p_n A_n + q_n J)^{-1} + r_n J \big)^{-1} \stackrel{H(\mathcal{M}(\Omega)^2)}{-\hspace{-0.2cm}\rightharpoonup}{p A_* + q J}.\end{equation}\label{lemDykhnelim}

\end{Prop}\begin{Rem}

Proposition \ref{lemDykhnelim} completes Theorem \ref{theo22} performed with the transformation \begin{equation}A \longmapsto \frac{A^\text{T}}{\det A} = J^{-1}A^{-1}J,\end{equation} to other Dykhne transformations of type (see \cite{MilBook}, Section 4.1): \begin{equation}
A \longmapsto \big((p A + q J)^{-1} + r J\big)^{-1} = \big(p A + q J\big)\big((1-rq) I_2 + rp JA\big)^{-1}
\end{equation}

\end{Rem}

\begin{Rem}

The convergence of $r_n$ to $r = 0$ is not necessary but sufficient for our purpose. If $r \neq 0$, the different convergences are conserved but lead us to the expression \begin{equation} p A_* + q J =\displaystyle B_* \big((1- qr) I_2 + p \ r J A_* \big). \end{equation}

\end{Rem}

\noindent \textbf{Proof of Proposition \ref{lemDykhnelim}}. The proof is divided into two steps. In the first step, we use Lemma \ref{divcurl} to show the $H(\mathcal{M}(\Omega)^2)$-convergence of $\widetilde{A}_n := p_n A_n + q_n J$ to $p A_* + q J$. In the second step, we build a matrix $Q_n$ which will be used as a corrector for $B_n$ and then use again Lemma \ref{divcurl}.

\bigskip

\noindent \textit{First step}: $\widetilde{A}_* = p A_* + q J$.

\smallskip

First of all, thanks to Theorem 2.2 \cite{NHM}, we already know that, up to a subsequence still denoted by $n$, $A_n$ $H(\mathcal{M}(\Omega)^2)$-converges to $A_*$. We consider a corrector $P_n$ associated with $A_n$ in the sense of Murat-Tartar (see, e.g., \cite{Mur2}), such that, for $\lambda \in \R^2$, $P_n \lambda = \nabla w_n^\lambda$ is defined by
\begin{equation}\; \left\{\!\! \begin{array}{l l}
\mathrm{div} (A_n \nabla w_n^\lambda ) = \mathrm{div} (A_* \nabla (\lambda \cdot x) )\vspace{0.2cm} \ \ & \mathrm{in} \ \Omega\\
w_n^\lambda = \lambda \cdot x \ \ & \mathrm{on} \ \partial \Omega\\
\end{array} \right. \;\;\;\end{equation}

\noindent Again with Theorem 2.2 of \cite{NHM} and Definition \hyperref[H(M)]{\ref*{H(M)}}, we know that $P_n \lambda$ converges weakly in $L^2(\Omega)^2$ to $\lambda$ and $A_n P_n \lambda$ converges weakly-$*$ in $\mathcal{M}(\Omega)$ to $A_* \lambda$.

Since, for any $\lambda, \mu \in \R^2$, \[ \alpha ||\nabla w_n^\mu||^2_{L^2(\Omega)^2} \leq \int_\Omega A_n \nabla w_n^\mu \cdot \nabla w_n^\mu \dx = \int_\Omega A_* \mu \cdot \nabla w_n^\mu \dx \leq 2 ||\bar{a}||_{L^\infty(\Omega)}|\mu| \ |\Omega|^{1/2} ||\nabla w_n^\mu||_{L^2(\Omega)^2}\] and\[\int_\Omega A_n^{-1} A_n \nabla w_n^\lambda \cdot A_n \nabla w_n^\lambda \dx = \int_\Omega A_n \nabla w_n^\lambda \cdot \nabla w_n^\lambda \dx,\]

\noindent the sequences $\xi_n := A_n \nabla w_n^\lambda$ and $v_n := w_n^\mu$ satisfy (\hyperref[divrot1]{\ref*{divrot1}}) and (\ref{divrot2}). This combined with (\ref{conddet}) implies that we can apply Lemma \ref{divcurl} to obtain\begin{equation}\forall\,\lambda, \mu \in \R, \quad A_n P_n \lambda \cdot P_n \mu \ \harpoon \widetilde{A}_* \lambda \cdot \mu \ \mathrm{in} \ \mathscr{D}'(\Omega).\label{convergence0}\end{equation}

We denote $\widetilde{A}_n := p_n A_n + q_n J$ and consider $\delta_n$ such that $\delta_n J:=A_n - A_n^s$. Then, the matrix $\widetilde{A}_n$ satisfies \begin{equation}\widetilde{A}_n \xi \cdot \xi = p_n A_n \xi \cdot \xi \geq p_n \alpha |\xi|^2 \label{equicoAtilde}.\end{equation} Moreover, \begin{equation*} \det \widetilde{A}_n= p_n^2 \det A_n^s + (p_n \delta_n + q_n)^2 \leq p_n^2 (\det A_n^s + 2 \delta_n^2) + 2 q_n^2 \leq 2 p_n^2 \det A_n + 2 q_n^2 \leq C \det A_n, \end{equation*}the last inequality being a consequence of $A_n \geq \alpha I_2$. This inequality gives, by \eqref{conddet}, \begin{equation}
\displaystyle{\frac{\det \widetilde{A}_n}{\det \widetilde{A}_n^s}} |\widetilde{A}_n^s| = \displaystyle{\frac{\det \widetilde{A}_n}{p_n^2 \det A_n^s}} p_n |A_n^s| \leq C \displaystyle{\frac{\det A_n}{ \det A_n^s}} |A_n^s| \leq C. \label{infdet}
\end{equation}

\noindent Then by \eqref{equicoAtilde}, \eqref{infdet} and \cite{NHM} again, up to a subsequence still denoted by $n$, $\widetilde{A}_n$ $H(\mathcal{M}(\Omega)^2)$-converges to $\widetilde{A}_*$ and we have, by the classical div-curl lemma of \cite{Mur2} for $J P_n \lambda \cdot P_n \mu$ and \eqref{convergence0}, \[ \forall\,\lambda, \mu \in \R, \quad (p_n A_n + q_n J) P_n \lambda \cdot P_n \mu = p_n A_n P_n \lambda \cdot P_n \mu + q_n J P_n \lambda \cdot P_n \mu \stackrel{\mathscr D'(\Omega)}{-\hspace{-0.2cm}\rightharpoonup} p A^* \lambda \cdot \mu + q J \lambda \cdot \mu,\]

\noindent that can be rewritten\[\widetilde{A}_* = p A_* + q J.\]

\bigskip

\noindent \textit{Second step}: $B_* = \widetilde{A}_*$.

\smallskip

Let $\theta \in \mathscr{C}_c^1(\Omega)$ and $\widetilde{P}_n$ a corrector associated with $\widetilde{A}_n$, such that, for $\lambda \in \R^2$, $\widetilde{P}_n \lambda = \nabla \widetilde{w}_n^\lambda$ is defined by
\begin{equation}\; \left\{\!\! \begin{array}{l l}
\mathrm{div} (\widetilde{A}_n \nabla \widetilde{w}_n^\lambda ) = \mathrm{div} \big(\widetilde{A}_* \nabla (\theta \lambda \cdot x)\big)\vspace{0.2cm} \ \ & \mathrm{in} \ \Omega\\
\widetilde{w}_n^\lambda = 0 \ \ & \mathrm{on} \ \partial \Omega.\\
\end{array} \right. \;\;\;\label{eqw_n}\end{equation}

\noindent By Definition \ref{H(M)}, we have \begin{equation}\; \left\{\!\! \begin{array}{ r l l}
\widetilde{w}_n^\lambda & \harpoon 0\theta \lambda \cdot x &\ \ \mathrm{weakly} \ \mathrm{in} \ H_0^1(\Omega)\vspace{0.2cm}, \\
\widetilde{A}_n \nabla \widetilde{w}_n^\lambda & \harpoon 0\widetilde{A}_* \nabla (\theta \lambda \cdot x) &\ \ \text{weakly-$*$} \ \mathrm{in} \ \mathcal M (\Omega)^2.\\
\end{array} \right. \;\;\;\label{convcorrtilde}\end{equation}

Let us now consider $B_n = \big( \widetilde{A}_n^{-1} + r_n J \big)^{-1}$. $B_n$ is symmetric and so is its inverse : \[B_n^{-1} = \widetilde{A}_n^{-1} + r_n J = (\widetilde{A}_n^{-1} + r_n J)^s = (\widetilde{A}_n^{-1})^s. \] We then have, by a little computation (like in \eqref{calceasy}) and \eqref{infdet}, \begin{align} \displaystyle{\frac{\det B_n}{\det B_n^s}} |B_n^s| = |B_n| = \left| \left( \big( \widetilde{A}_n^{-1}\big)^s\right)^{-1} \right| = \frac{\det \widetilde{A}_n}{\det \widetilde{A}_n^s} |\widetilde{A}_n^s| \leq C. \label{Bn2}
\end{align}For any $\xi \in \R^2$, the sequence $\nu_n := (I + r_n J \widetilde{A}_n )^{-1} \xi$ satisfies, by \eqref{conddet}, \[ |\xi| \leq \left( 1 + ||r_n \widetilde{A}_n||_{L^\infty(\Omega)^{2 \times 2}} \right) |\nu_n| \leq \left( 1 + p_n ||r_n A_n||_{L^\infty(\Omega)^{2 \times 2}} + q_n r_n\right) |\nu_n| \leq (1+C) |\nu_n|, \] hence \begin{equation} B_n \xi \cdot \xi = \widetilde{A}_n \nu_n \cdot (I + r_n J \widetilde{A}_n ) \nu_n = \widetilde{A}_n \nu_n \cdot \nu_n = p_n A_n \nu_n \cdot \nu_n \geq p_n \alpha |\nu_n|^2 \geq \alpha \frac{p_n}{(1+C)^2} |\xi|^2 \geq C |\xi|^2 \label{Bn1}\end{equation}with $C >0$. Therefore, with \eqref{Bn1} and \eqref{Bn2}, again by Theorem 2.2 of \cite{NHM}, up to a subsequence still denoted by $n$, $B_n$ $H(\mathcal{M}(\Omega)^2)$-converges to $B_*$.

Let $\psi \in \mathscr C_c^1(\Omega)$ and $R_n$ be a corrector associated to $B_n$, such that, for $\mu \in \R^2$, $R_n \mu = \nabla v_n^\mu$ is defined by
\begin{equation}\; \left\{\!\! \begin{array}{ l l}
\mathrm{div} \left(B_n \nabla v_n^\mu \right) = \mathrm{div} \big(B_* \nabla (\psi \mu \cdot x)\big)\vspace{0.2cm} &\ \ \mathrm{in} \ \Omega\\
v_n^\mu = 0 &\ \ \mathrm{on} \ \partial \Omega.\\
\end{array} \right. \;\;\; \label{Rn}
\end{equation} By Definition \ref{H(M)}, we have the convergences\begin{equation}\; \left\{\!\! \begin{array}{ r l l}
v_n^\mu & \harpoon 0\psi \mu \cdot x \vspace{0.2cm} &\ \ \mathrm{weakly} \ \mathrm{in} \ H_0^1(\Omega), \\
B_n \nabla v_n^\mu & \harpoon 0 B_* \nabla (\psi \mu \cdot x) &\ \ \text{weakly-$*$} \ \mathrm{in} \ \mathcal M (\Omega)^2.\\
\end{array} \right. \;\;\; \label{convvn}\end{equation}Let us define the matrix $Q_n := (I + r_n J \widetilde{A}_n) \widetilde{P}_n$. We have\begin{align}
B_n Q_n = (\widetilde{A}_n^{-1} + r_n J)^{-1} (I + r_n J \widetilde{A}_n) \widetilde{P}_n = (\widetilde{A}_n^{-1} + r_n J)^{-1} (\widetilde{A}_n^{-1} + r_n J) \widetilde{A}_n \widetilde{P}_n = \widetilde{A}_n \widetilde{P}_n \label{BQ = AP}.\end{align}

We are going to pass to the limit in $\mathscr D'(\Omega)$ the equality given by \eqref{BQ = AP} and the symmetry of $B_n$:\begin{equation}
\widetilde{A}_n \widetilde{P}_n \lambda \cdot R_n \mu = B_n Q_n \lambda \cdot R_n \mu = Q_n \lambda \cdot B_n R_n \mu \label{egall}.
\end{equation}

\noindent On the one hand, $\widetilde{A}_n$ satisfies \eqref{conddet1} by \eqref{equicoAtilde} and \eqref{infdet}. The sequences $\xi_n := \widetilde{A}_n \widetilde{P}_n \lambda$ and $v_n :=v_n^\mu$ satisfy the hypothesis \eqref{divrot2} by \eqref{eqw_n} and \eqref{divrot1} because \begin{equation*}
\int_\Omega \big(\widetilde{A}_n\big)^{-1} \xi_n \cdot \xi_n \dx + ||v_n||_{H_0^1(\Omega)} = \int_\Omega \widetilde{A}_n \widetilde{P}_n \lambda \cdot \widetilde{P}_n \lambda \dx + ||v_n^\mu||_{H_0^1(\Omega)} \dx \leq C
\end{equation*} by \eqref{eqw_n} and the convergences \eqref{convvn} and \eqref{convcorrtilde}. The application of Lemma \ref{divcurl}, \eqref{convcorrtilde} and \eqref{convvn} give the convergence\begin{equation}\widetilde{A}_n \widetilde{P}_n \lambda \cdot R_n \mu \ \harpoon \widetilde{A}^* \nabla (\theta \lambda \cdot x) \cdot \nabla (\psi \mu \cdot x) \quad \mathrm{in} \ \mathscr{D}'(\Omega). \label{convA*}\end{equation}

\noindent On the other hand, we have the equality\begin{equation}
Q_n \lambda \cdot B_n R_n \mu = B_n R_n \mu \cdot \widetilde{P}_n \lambda + B_n R_n \mu \cdot r_n J \widetilde{A}_n \widetilde{P}_n \label{egalBn}.
\end{equation}The matrix $B_n$ satisfies \eqref{conddet1} by \eqref{Bn1} and \eqref{Bn2}. The sequences $\xi_n := B_n R_n \mu$ and $v_n :=\widetilde{w}_n^\lambda$ satisfy the hypothesis \eqref{divrot2} by \eqref{Rn} and \eqref{divrot1} of Lemma \ref{divcurl} because \begin{equation*}
\int_\Omega \big(B_n\big)^{-1} \xi_n \cdot \xi_n \dx + ||v_n||_{H_0^1(\Omega)} = \int_\Omega B_n R_n \mu \cdot R_n \mu \dx + ||\widetilde{w}_n^\lambda||_{H_0^1(\Omega)} \dx \leq C
\end{equation*} by \eqref{Rn} and the convergences \eqref{convcorrtilde} and \eqref{convvn}. The application of Lemma \ref{divcurl}, \eqref{convcorrtilde} and \eqref{convvn} give the convergence\begin{equation}B_n R_n \mu \cdot \widetilde{P}_n \lambda \ \harpoon 0 B_* \nabla (\psi \mu \cdot x) \cdot \nabla (\theta \lambda \cdot x) \quad \mathrm{in} \ \mathscr{D}'(\Omega). \label{convB*1}\end{equation}The convergence of the right part of \eqref{egalBn} is more delicate. The demonstration is the same as for Lemma \ref{divcurl}. Let $\omega$ be a simply connected open subset of $\Omega$ such as $\omega \subset \subset \Omega$. The function $ \widetilde{A}_n \widetilde{P}_n \lambda - \widetilde{A}_* \nabla (\theta \lambda \cdot x)$ is divergence-free and we can introduce a function $z_n^\lambda$ such as \begin{align}\widetilde{A}_n \widetilde{P}_n \lambda = \widetilde{A}_* \nabla (\theta \lambda \cdot x) + J \nabla z_n^\lambda, \\ z_n^\lambda \longrightarrow 0 \quad \mathrm {strongly} \ \mathrm{in} \ L_{\textrm{loc}}^2(\omega). \label{zlambda}\end{align}The equality \begin{align*} B_n R_n \mu \cdot r_n J \widetilde{A}_n \widetilde{P}_n \lambda &= r_n B_n R_n \mu \cdot J \widetilde{A}_* \nabla (\theta \lambda \cdot x) - r_n B_n R_n \mu \cdot \nabla z_n^\lambda \\
&= r_n B_n R_n \mu \cdot J \widetilde{A}_* \nabla (\theta \lambda \cdot x) - r_n \textrm{div} ( z_n^\lambda B_n R_n \mu) + r_n z_n^\lambda \ \mathrm{div}\big( B_* \nabla (\theta \lambda \cdot x) \big) \end{align*} leads us, by \eqref{convvn}, \eqref{zlambda} and the convergence to $0$ of $r_n$, like in the demonstration of Lemma \ref{divcurl}, to\begin{equation}B_n R_n \mu \cdot r_n J \widetilde{A}_n \widetilde{P}_n \ \harpoon 00 \quad \mathrm{in} \ \mathscr{D}'(\omega). \label{convB*2}\end{equation}

Finally, by combining \eqref{egall}, \eqref{convA*}, \eqref{convB*1} and \eqref{convB*2}, we obtain, for any simply connected open subset $\omega$ of $\Omega$ such as $\omega \subset \subset \Omega$, \[\widetilde{A}_* \nabla (\theta \lambda \cdot x) \cdot \nabla (\psi \mu \cdot x) = B_* \nabla (\psi \mu \cdot x) \cdot \nabla (\theta \lambda \cdot x) \quad \mathrm{in} \ \mathscr{D}'(\omega).\] We conclude, by taking $\theta = 1$ and $\psi =1$ on $\omega$ and taking into account that $B_*$ is symmetric and $\omega$, $\lambda$, $\mu$ are arbitrary, that:\[B_* = \widetilde{A}_* = p A_* + qJ.\] \qed

\subsection{An application to isotropic two-phase media}

In this section, we study the homogenization of a two-phase isotropic medium with high contrast and non-necessarily symmetric conductivities. The study of the symmetric case in Proposition \ref{diag} permits to obtain Theorem \ref{divcurl1} by applying the transformation of Proposition \ref{lemDykhnelim}. We use Notation \ref{defabs}.

\begin{Prop}

Let $\Omega$ be a bounded open subset of $\R^2$ such that $|\partial \Omega| = 0 $. Let $\omega_n$, $n$ in $\mathbb N$, be a sequence of open subsets of $\Omega$ with characteristic function $\chi_n$, satisfying $\theta_n:= |\omega_n| < 1$, $\theta_n$ converges to $0$, and \begin{equation}
\displaystyle \frac{\chi_n}{\theta_n} \xrightharpoonup{\textcolor{white}{H(}} a \in L^\infty(\Omega) \quad \text{weakly-$*$ in} \ \mathcal M(\Omega)\label{conda}.\end{equation}We assume that there exists $\alpha_1, \alpha_2 > 0$ and two positive sequences $\alpha_{1,n}$, $\alpha_{2,n} \geq a_0 >0$ verifying \begin{equation} \lim \limits_{n \to \infty} \alpha_{1,n} = \alpha_1 \quad \text{and} \quad  \lim \limits_{n \to \infty}\theta_n\alpha_{2,n} = \alpha_2, \label{cond88}\end{equation} and that the conductivity takes the form $$\sigma_n^0(\alpha_{1,n},\alpha_{2,n}) = (1-\chi_n) \alpha_{1,n}I_2 + \chi_n \alpha_{2,n}I_2.$$
\noindent Then, there exists a subsequence of $n$, still denoted by $n$, and a locally Lipschitz function \[\sigma_*^0 : (0 , \infty)^2 \longrightarrow \mathcal M(a_0, 2 || a ||_{\infty};\Omega)\] such that \begin{eqnarray}\forall\,(\alpha_1,\alpha_2) \in (0,\infty)^2, \quad \sigma_n^0(\alpha_{1,n},\alpha_{2,n}) \stackrel{H(\mathcal{M}(\Omega)^2)}{-\hspace{-0.2cm}\rightharpoonup}{\sigma_*^0(\alpha_1,\alpha_2)}\label{defsigma*}.\end{eqnarray}
\label{diag}
\end{Prop}

\noindent \textbf{Proof of Proposition \ref{diag}}. The proof is divided into two parts. We first prove the theorem for $\alpha_{1,n} = \alpha_1$, $\alpha_{2,n} = \theta_n^{-1} \alpha_2$, and then treat the general case. 

\bigskip

\noindent \textit{First step}: The case $\alpha_{1,n} = \alpha_1$, $\alpha_{2,n} = \theta_n^{-1} \alpha_2$.
\par\noindent \smallskip
In this step we denote $\sigma_n^0(\alpha) := \sigma_n^0(\alpha_{1},\theta_n^{-1}\alpha_{2})$, for $\alpha = (\alpha_1,\alpha_2) \in (0,\infty)^2$.
Theorem 2.2 of \cite{NHM} implies that for any $\alpha \in (0,\infty)^2$, there exists a subsequence of $n$ such that $\sigma^0_n(\alpha)$ $H(\mathcal{M}(\Omega)^2)$-converges in the sense of Definition \ref{H(M)} to some matrix-valued function in $\mathcal M(a_0, 2 || a ||_{\infty};\Omega)$.

By a diagonal extraction, there exists a subsequence of $n$, still denoted by $n$, such that\begin{equation}\forall\,\alpha \in \mathbb Q^2 \cap (0, \infty)^2, \quad \sigma^0_n(\alpha)\stackrel{H(\mathcal{M}(\Omega)^2)}{-\hspace{-0.2cm}\rightharpoonup}{\sigma^0_*(\alpha)}. \label{convbase}\end{equation}We are going to show that this convergence is true any pair $\alpha \in (0,\infty)^2$.

We have, by \eqref{conda}, for any $\alpha \in \mathbb Q^2 \cap (0,\infty)^2$, \begin{equation}\displaystyle |\sigma^0_n(\alpha)| = (1-\chi_n)\alpha_1 + \chi_n \displaystyle{\frac{\alpha_2}{\theta_n}} \xrightharpoonup{\textcolor{white}{H(}} \alpha_1 + \alpha_2 \ a \in L^\infty(\Omega) \quad \text{ weakly-$*$ in} \ \mathcal M(\Omega)\label{aa1}\end{equation} and, since $\theta_n \in (0,1)$, \begin{equation} \forall\,\xi \in \R^2, \quad \sigma^0_n(\alpha) \xi \cdot\xi = \alpha_1 (1-\chi_n) |\xi|^2 + \chi_n \displaystyle{\frac{\alpha_2}{\theta_n}}|\xi|^2 \geq \min(\alpha_1,\alpha_2) |\xi|^2 \quad \text{a.e. in} \ \Omega. \label{inegamma} \end{equation}By applying Theorem 2.2 of \cite{NHM} with \eqref{aa1}, we have the inequality\begin{equation} |\sigma^0_*(\alpha) \lambda| \leq 2 |\lambda| \left(\alpha_1 + \alpha_2 ||a||_\infty \right).\label{ineqa}\end{equation}

For any $\alpha \in \mathbb Q^2 \cap (0,\infty)^2$ and $\lambda \in \R^2$, consider the corrector $w_n^{\alpha,\lambda}$ associated with $\sigma^0_n(\alpha)$ defined by
\begin{equation}\; \left\{\!\! \begin{array}{ r c l l}
\mathrm{div} \left(\sigma^0_n(\alpha) \nabla w_n^{\alpha,\lambda} \right) &=& \mathrm{div} \left(\sigma^0_*(\alpha) \lambda \right) \vspace{0.2cm}&\ \ \mathrm{in} \ \Omega,\\
w_n^{\alpha,\lambda} &=& \lambda \cdot x &\ \ \mathrm{on} \ \partial \Omega,\\
\end{array} \right. \;\;\;\label{corra}\end{equation}which depends linearly on $\lambda$.

Let $\alpha \in \mathbb Q^2 \cap (0,\infty)^2$. Let us show that the energies \begin{equation}\int_\Omega \sigma^0_n(\alpha) \nabla w_n^{\alpha,\lambda} \cdot \nabla w_n^{\alpha,\lambda} \dx\label{ener11}\end{equation}are bounded. We have, by \eqref{corra}, \eqref{ineqa} and the Cauchy-Schwarz inequality\begin{align*}
& \int_\Omega \sigma^0_n(\alpha) \nabla w_n^{\alpha,\lambda} \cdot \nabla w_n^{\alpha,\lambda} \dx \\
& = \int_\Omega \sigma^0_*(\alpha) \lambda \cdot \big(\nabla w_n^{\alpha,\lambda} - \lambda \big) \dx + \int_\Omega \sigma^0_n(\alpha) \nabla w_n^{\alpha,\lambda} \cdot \lambda \dx\\
& = \int_\Omega \sigma^0_*(\alpha) \lambda \cdot \nabla w_n^{\alpha,\lambda}\dx -\int_\Omega \underbrace{\sigma^0_*(\alpha) \lambda \cdot \lambda}_{\geq 0} \dx + \int_\Omega \sigma^0_n(\alpha) \nabla w_n^{\alpha,\lambda} \cdot \lambda \dx
\end{align*}which leads us to \begin{equation}
\int_\Omega \sigma^0_n(\alpha) \nabla w_n^{\alpha,\lambda} \cdot \nabla w_n^{\alpha,\lambda} \dx \leq \int_\Omega |\sigma^0_*(\alpha) \lambda \cdot \nabla w_n^{\alpha,\lambda}|\dx + \int_\Omega |\sigma^0_n(\alpha) \nabla w_n^{\alpha,\lambda} \cdot \lambda| \dx. \label{ener}
\end{equation}

On the one hand, the Cauchy-Schwarz inequality gives \begin{equation*}
\left(\int_\Omega |\sigma^0_n(\alpha) \nabla w_n^{\alpha,\lambda} \cdot \lambda| \dx \right)^2 \leq |\lambda|^2 \int_\Omega |\sigma^0_n(\alpha)| \dx \int_\Omega \sigma^0_n(\alpha) \nabla w_n^{\alpha,\lambda} \cdot \nabla w_n^{\alpha,\lambda} \dx
\end{equation*} that is\begin{equation}\left(\int_\Omega |\sigma^0_n(\alpha) \nabla w_n^{\alpha,\lambda} \cdot \lambda| \dx \right)^2 \leq |\lambda|^2 |\alpha| \int_\Omega \sigma^0_n(\alpha) \nabla w_n^{\alpha,\lambda} \cdot \nabla w_n^{\alpha,\lambda} \dx. \label{001}
\end{equation}

On the other hand, by \eqref{inegamma} and the Cauchy-Schwarz inequality, we have\begin{align*} \int_\Omega |\sigma^0_*(\alpha) \lambda \cdot \nabla w_n^{\alpha,\lambda}|\dx &\leq 2 |\lambda| (\alpha_1 + \alpha_2 ||a||_\infty) \sqrt{\int_\Omega |\nabla w_n^{\alpha,\lambda}|^2 \dx} \\
&\leq 2 |\lambda| (\alpha_1 + \alpha_2 ||a||_\infty) \sqrt{\frac{1}{\alpha_1} + \frac{1}{\alpha_2}} \sqrt{ \int_\Omega \sigma^0_n(\alpha) \nabla w_n^{\alpha,\lambda} \cdot \nabla w_n^{\alpha,\lambda} \dx} \end{align*}that is\begin{equation} \int_\Omega |\sigma^0_*(\alpha) \lambda \cdot \nabla w_n^{\alpha,\lambda}|\dx \leq C \ |\lambda|^2 \ |\alpha| \sqrt{\frac{1}{\alpha_1} + \frac{1}{\alpha_2}} \sqrt{ \int_\Omega \sigma^0_n(\alpha) \nabla w_n^{\alpha,\lambda} \cdot \nabla w_n^{\alpha,\lambda} \dx} \label{ineqlips}\end{equation} where $C$ does not depend on $n$ nor $\alpha$.

By combining \eqref{ener}, \eqref{001} and \eqref{ineqlips}, we have \begin{equation}
\int_\Omega \sigma^0_n(\alpha) \nabla w_n^{\alpha,\lambda} \cdot \nabla w_n^{\alpha,\lambda} \dx \leq C \ |\lambda|^2 \ \big(\underbrace{ |\alpha| + |\alpha|^2 ({\alpha_1}^{-1} + {\alpha_2}^{-1})}_{=: M(\alpha)} \big) \label{ineqener}
\end{equation}where $C$ does not depend on $n$ nor $\alpha$.

Let ${\alpha}' \in \mathbb Q^2 \cap (0,\infty)^2$. The sequences $\xi_n := \sigma^0_n(\alpha) \nabla w_n^{\alpha,\lambda}$ and $v_n : = w_n^{{\alpha}',\lambda}$ satisfy the assumptions \eqref{divrot1} and \eqref{divrot2} of Lemma \ref{divcurl}. By symmetry, we have the convergences \begin{equation}\; \left\{\!\! \begin{array}{ r l l}
\sigma^0_n(\alpha) \nabla w_n^{\alpha,\lambda} \cdot \nabla w_n^{{\alpha}',\lambda}& \harpoon 0\sigma^0_*(\alpha) \lambda \cdot \lambda \vspace{0.2cm} & \ \mathrm{weakly} \ \mathrm{in} \ \mathscr D'(\Omega), \\
\sigma^0_n({\alpha}') \nabla w_n^{{\alpha}',\lambda} \cdot \nabla w_n^{\alpha,\lambda}& \harpoon 0\sigma^0_*({\alpha}') \lambda \cdot \lambda \vspace{0.2cm} & \ \mathrm{weakly} \ \mathrm{in} \ \mathscr D'(\Omega).
\end{array} \right. \;\;\; \label{crossconv}\end{equation} As the matrices are symmetric, we have\begin{equation*}
\left(\sigma^0_n(\alpha)-\sigma^0_n({\alpha}')\right) \nabla w_n^{\alpha,\lambda} \cdot \nabla w_n^{{\alpha}',\lambda} = \sigma^0_n(\alpha) \nabla w_n^{\alpha,\lambda} \cdot \nabla w_n^{{\alpha}',\lambda} - \sigma^0_n({\alpha}') \nabla w_n^{{\alpha}',\lambda} \cdot \nabla w_n^{\alpha,\lambda},
\end{equation*}hence\begin{equation}
\left(\sigma^0_n(\alpha)-\sigma^0_n({\alpha}')\right) \nabla w_n^{\alpha,\lambda} \cdot \nabla w_n^{{\alpha}',\lambda} \ \harpoon 0\left(\sigma^0_*(\alpha)-\sigma^0_*({\alpha}')\right) \lambda \cdot \lambda \quad \mathrm{weakly} \ \mathrm{in} \ \mathscr D'(\Omega). \label{convdiffs}
\end{equation}

Let $\lambda \in \R^2$. We have, by the Cauchy-Schwarz inequality, with the Einstein convention \begin{align*}
& \int_\Omega \big|\big(\sigma^0_n(\alpha) - \sigma^0_n(\alpha') \big) \nabla w_n^{\alpha,\lambda} \cdot \nabla w_n^{\alpha',\lambda}\big| \dx \\
 & = \int_{\Omega \setminus \omega_n} |\alpha_1 - {\alpha}'_1| \ \big| \nabla w_n^{{\alpha},\lambda} \cdot \nabla w_n^{{\alpha}',\lambda} \big| \dx + \fint_{\omega_n} |\alpha_2 - {\alpha}'_2| \ \big| \nabla w_n^{{\alpha},\lambda} \cdot \nabla w_n^{{\alpha}',\lambda} \big| \dx \\
 & \leq  |\alpha_1 - {\alpha}'_1| \sqrt{\int_{\Omega\setminus\omega_n} |\nabla w_n^{{\alpha},\lambda}|^2 \dx } \ \sqrt{\int_{\Omega\setminus\omega_n} |\nabla w_n^{{{\alpha}'},\lambda}|^2 \dx } \\
& \hspace{0.5cm} +  |\alpha_2 - {\alpha}'_2| \sqrt{\fint_{\omega_n} |\nabla w_n^{{\alpha},\lambda}|^2 \dx } \ \sqrt{\fint_{\omega_n} |\nabla w_n^{{{\alpha}'},\lambda}|^2 \dx } \\
&\leq |\alpha_i - {\alpha}'_i| \sqrt{\frac{1}{\alpha_i} \int_\Omega \sigma^0_n(\alpha) \nabla w_n^{\alpha,\lambda} \cdot \nabla w_n^{\alpha,\lambda} \dx } \ \sqrt{\frac{1}{{\alpha}'_i} \int_\Omega \sigma^0_n(\alpha) \nabla w_n^{{\alpha}',\lambda} \cdot \nabla w_n^{{\alpha}',\lambda} \dx } .
\end{align*} This combined with \eqref{ineqener} yields\begin{equation*}
\int_\Omega \big|\big(\sigma^0_n(\alpha) - \sigma^0_n(\alpha')\big) \nabla w_n^{\alpha,\lambda} \cdot \nabla w_n^{\alpha',\lambda}\big| \leq C |\lambda|^2 \ \frac{|\alpha_i - {\alpha}'_i|}{\sqrt{|\alpha_i| |{\alpha}'_i| } } \ M(\alpha) \ M({\alpha}')
\end{equation*} The sequence of \eqref{convdiffs} is thus bounded in $L^1(\Omega)^2$ which implies that \eqref{convdiffs} holds weakly-$*$ in $\mathcal M(\Omega)$. Hence, we get, for any $\varphi \in \mathscr C_c(\Omega)$, that\begin{equation}
\int_\Omega \left|\left(\sigma^0_*(\alpha) - \sigma^0_*({\alpha}') \right) \lambda \cdot \lambda \right| \ \varphi \dx  \leq C \ |\lambda|^2 \ \frac{|\alpha_i - {\alpha}'_i|}{\sqrt{|\alpha_i| |{\alpha}'_i| } } \ M(\alpha) \ M({\alpha}')  \ || \varphi||_\infty. \label{ineqsigma}
\end{equation}

\noindent Then, the Riesz representation theorem implies that \begin{equation*} \left| \left| \sigma^0_*(\alpha) - \sigma^0_*({\alpha}') \right| \right|_{L^1(\Omega)^{2 \times 2}} \leq C \ \frac{|\alpha_i - {\alpha}'_i|}{\sqrt{|\alpha_i| |{\alpha}'_i| } } \ M(\alpha) \ M(\alpha').\end{equation*}Therefore, by the definition of $M$ in \eqref{ineqener}, for any compact subset $K \subset (0,\infty)^2$,\begin{equation} \exists\,C > 0, \quad \forall\,\alpha, {\alpha}' \in \mathbb Q^2 \cap K, \quad \left| \left| \sigma^0_*(\alpha) - \sigma^0_*({\alpha}') \right| \right|_{L^1(\Omega)^{2 \times 2}} \leq C \ |\alpha - {\alpha}'| . \label{ineqlips1} \end{equation} This estimate permits to extend the definition \eqref{convbase} of $\sigma_*^0$ on $(0,\infty)^2$ by \begin{equation} \forall\,\alpha \in (0, \infty)^2, \quad\sigma_*^0(\alpha)  = \lim_{\substack{\alpha'\to \alpha \\ \alpha' \in \mathbb Q^2 \cap (0, \infty)^2}} \sigma_*^0(\alpha') \quad \text{strongly in} \ L^1(\Omega)^{2 \times 2}. \label{defext}\end{equation}

Let $\alpha \in (0, \infty)^2$. Theorem 2.2 of \cite{NHM} implies that there exists a subsequence of $n$, denoted by $n'$, and a matrix-valued function $\widetilde{\sigma}_* \in \mathcal M(a_0, 2 || a ||_{\infty};\Omega)$ such that \begin{equation} \sigma_{n'}(\alpha) \stackrel{H(\mathcal{M}(\Omega)^2)}{-\hspace{-0.2cm}\rightharpoonup}{\widetilde{\sigma}_*} \label{uu} .\end{equation} Repeating the arguments leading to \eqref{ineqlips1}, for any positive sequence of rational pair $(\alpha^q)_{q \in \mathbb N}$ converging to $\alpha$, we have\begin{equation} \exists\,C > 0, \quad \left| \left| \widetilde{\sigma}_* - \sigma^0_*({\alpha^q}) \right| \right|_{L^1(\Omega)^{2 \times 2}} \leq C \ |\alpha - \alpha^q| , \label{ineqlips2} \end{equation} hence, by \eqref{defext}, $\widetilde{\sigma}_* = \sigma^0_*(\alpha)$. Therefore by the uniqueness of the limit in \eqref{uu}, we obtain for the whole sequence satisfying \eqref{convbase} \begin{equation}\forall\,\alpha \in (0,\infty)^2, \quad \sigma_{n}(\alpha) \stackrel{H(\mathcal{M}(\Omega)^2)}{-\hspace{-0.2cm}\rightharpoonup}{\sigma^0_*(\alpha)}.\end{equation}In particular, the function $\sigma^0_*$ satisfies \eqref{ineqlips1} and \eqref{defext}, i.e. $\sigma^0_*$ is a locally Lipschitz function on $(0,\infty)^2$.

\bigskip

\noindent \textit{Second step}: The general case.
\par\noindent \smallskip
We denote $\alpha^n = (\alpha_{1,n},\alpha_{2,n})$ and $ \sigma_n^0(\alpha^n) = \sigma_n^0(\alpha_{1,n},\alpha_{2,n})$. Theorem 2.2 of \cite{NHM} implies that there exists a subsequence of $n$, denoted by $n'$, such that $\sigma^0_{n'}(\alpha^{n'})$ $H(\mathcal{M}(\Omega)^2)$-converges to some \textcolor{white}{tt} ${\widetilde{\sigma}_* \in \mathcal M(a_0, 2 || a ||_{\infty};\Omega)}$ in the sense of Definition \ref{H(M)}.

As in the first step, for any $\alpha^{n'} \in (0,\infty)^2$ and $\lambda \in \R^2$, we can consider the corrector $w_{n'}^{\alpha^{n'},\lambda}$ associated with $\sigma^0_{n'}(\alpha^{n'})$ defined by
\begin{equation}\; \left\{\!\! \begin{array}{ r c l l}
\mathrm{div} \left(\sigma^0_{n'}(\alpha^{n'}) \nabla w_{n'}^{\alpha^{n'},\lambda} \right) &=& \mathrm{div} \left(\widetilde{\sigma}_* \lambda \right) \vspace{0.2cm}&\ \ \mathrm{in} \ \Omega,\\
w_{n'}^{\alpha^{n'},\lambda} &=& \lambda \cdot x &\ \ \mathrm{on} \ \partial \Omega,\\
\end{array} \right. \;\;\;\end{equation}which depends linearly on $\lambda$. Proceeding as in the first step, we obtain like in \eqref{convdiffs}, with $\alpha = (\alpha_1, \alpha_2)$ the limit of $\alpha^n$ according to \eqref{cond88},\begin{equation}
\left(\sigma^0_{n'}(\alpha)-\sigma^0_{n'}({\alpha^{n'}})\right) \nabla w_{n'}^{\alpha^{n'},\lambda} \cdot \nabla w_{n'}^{{\alpha},\lambda} \ \harpoon 0\left(\sigma^0_*(\alpha)-\widetilde{\sigma}_*\right) \lambda \cdot \lambda \quad \mathrm{weakly} \ \mathrm{in} \ \mathscr D'(\Omega). \label{convdiffs2}
\end{equation} Moreover, by the energy bound \eqref{ineqener}, which also holds for $\alpha^{n'}$, we have, for any $\varphi \in \mathscr D(\Omega)$, \begin{equation*}
\int_\Omega \left(\sigma^0_{n'}(\alpha)-\sigma^0_{n'}({\alpha^{n'}})\right) \nabla w_{n'}^{\alpha^{n'},\lambda} \cdot \nabla w_{n'}^{{\alpha},\lambda} \ \varphi \dx \underset{{n'} \to \infty}{\longrightarrow} 0.
\end{equation*}

This combined with \eqref{convdiffs2}, yields \begin{equation*}
\int_\Omega \left(\sigma^0_*(\alpha) - \widetilde{\sigma}_*\right) \lambda \cdot \lambda \ \varphi \dx = 0,
\end{equation*} which implies that $\sigma^0_*(\alpha) = \widetilde{\sigma}_*$. We conclude by a uniqueness argument. \qed

\bigskip

We can now obtain a result for (perturbed) non-symmetric conductivities. Then, we will use a Dykhne transformation to recover the symmetric case following the Milton approach \cite{MilBook} (pp. 61--65). This will allow us to apply Proposition \ref{diag}.

\begin{Theo}

Let $\Omega$ be a bounded open subset of $\R^2$ such that $|\partial \Omega| = 0 $. Let $\omega_n$, $n \in \mathbb N$, be a sequence of open subsets of $\Omega$ and denote by $\chi_n$ their characteristic function. We assume that $\theta_n = |\omega_n| < 1$ converges to $0$ and\begin{equation}
\displaystyle \frac{\chi_n}{\theta_n} \xrightharpoonup{\textcolor{white}{H(}} a \in L^\infty(\Omega) \quad \text{weakly-$*$ in} \ \mathcal M(\Omega)\label{conda'}.\end{equation}Consider the conductivity defined by\begin{eqnarray}\sigma_n(h) = (1-\chi_n) \sigma_1(h) + \displaystyle{\frac{\chi_n}{\theta_n}} \sigma_2(h)\end{eqnarray}

\noindent where for $j = 1,2$, $\sigma_j(h) = \alpha_j + h \beta_j J \in \mathbb R^{2 \times 2}$ with $\alpha_1,\alpha_2 > 0$ and $(\beta_1, \beta_2) \neq (0,0)$.

\noindent Then, there exists a subsequence of $n$, still denoted by $n$, and a locally Lipschitz function \[\sigma^0_* : (0 , \infty)^2 \longrightarrow \mathcal M \Big(\min( \alpha_1 , \alpha_2) , 2 \big( |\sigma_1| + |\sigma_2| \ ||a||_\infty \big) ; \Omega \Big)\]such that$$\sigma_n(h) \stackrel{H(\mathcal{M}(\Omega)^2)}{-\hspace{-0.2cm}\rightharpoonup}{\displaystyle \sigma^0_*\big(\alpha_1, \alpha_2 + \alpha_2^{-1} \beta_2^2 h^2 \big)} + h \beta_1 J.$$
\label{divcurl1}
\end{Theo}

\noindent \textbf{Proof of Theorem \ref{divcurl1}}. We have \[ \forall\,\xi \in \R^2, \quad \sigma_n(h) \xi \cdot \xi = (1-\chi_n) \alpha_1 |\xi|^2 + \displaystyle{\frac{\chi_n}{\theta_n}} \alpha_2 |\xi|^2 \geq \min(\alpha_1 , \alpha_2) |\xi|^2 \quad \text{a.e. in} \ \Omega\]and, by \eqref{conda'}, \[|\sigma_n(h)| = (1-\chi_n) |\sigma_1(h)| + \displaystyle{\frac{\chi_n}{\theta_n}} |\sigma_2(h)| \harpoon 0 |\sigma_1(h)| + a |\sigma_2(h)| \in L^\infty(\Omega) \quad \text{weakly-$*$ in} \ \mathcal M(\Omega).\]

In order to make a Dykhne transformation like in p.62 of \cite{MilBook}, we consider two real coefficients $a_n$ and $b_n$ in such a way that \[B_n := \big( a_n \sigma_n(h) + b_n J \big) \big( a_n I_2 + J \sigma_n(h) \big)^{-1} =  \left( ( p_n \sigma_n(h) + q_n J)^{-1} + r_n J \right)^{-1}\] is symmetric. An easy computation shows that the previous equality holds when \[ \displaystyle p_n := \frac{a_n^2}{a_n^2 + b_n}, \quad q_n := \frac{a_n b_n}{a_n^2 + b_n} \quad \textrm{and} \quad r_n :=  \frac{1}{a_n}.\] On the one hand, the estimates \eqref{g3} and \eqref{g22} with $\alpha_{2,n} = \theta_n^{-1} \alpha_2$, $\beta_{2,n} = \theta_n^{-1} \beta_2$, yield (note that they are independent of $\chi_n$) \begin{equation} \displaystyle p_n \underset{n \to \infty}{\sim} 1, \quad q_n \underset{n \to \infty}{\longrightarrow} - h \beta_1, \quad  \ r_n \underset{n \to \infty}{\longrightarrow} 0 \quad \text{and} \quad ||r_n \sigma_n(h)||_\infty \leq C \big(|\sigma_1(h)| + |\sigma_2(h)| \big). \label{4444}\end{equation}

On the other hand, as in Section \ref{2phasres}, with Notation \ref{defabs} and \eqref{zero}, we have\begin{equation}
B_n =  \sigma_n^0\big(\alpha'_{1,n}(h),\alpha'_{2,n}(h)\big), \label{zero1}
\end{equation} where \begin{equation} \alpha'_{1,n}(h) = \frac{a_n(\alpha_{1} + i h \beta_1) + i b_n}{a_n + i(\alpha_{1} + i h \beta_1)} \quad \textrm{and} \quad \alpha'_{2,n}(h) = \frac{a_n(\alpha_{2}/\theta_n + i h \beta_{2}/\theta_n) + i b_n}{a_n + i(\alpha_{2}/\theta_n + i h \beta_{2}/\theta_n)}.\end{equation}Hence, like in \eqref{2phas}, we have \begin{equation} \lim \limits_{n \to \infty} \alpha'_{1,n}(h) =\alpha_1 \quad \text{and} \quad  \lim \limits_{n \to \infty} \theta_n\alpha'_{2,n}(h) = \alpha_2 + \alpha_2^{-1} {\beta_2}^2 h^2. \label{compphas1}\end{equation}

We can first apply Proposition \ref{diag} with the conditions \eqref{zero1} and \eqref{compphas1} to have the $H( \mathcal M({\Omega})^2)$-convergence of $B_n$. Then, by virtue of Proposition \ref{lemDykhnelim}, with \eqref{4444} we get that $$\sigma_n(h) \stackrel{H(\mathcal{M}(\Omega)^2)}{-\hspace{-0.2cm}\rightharpoonup}{\displaystyle{\sigma^0_*\big(\alpha_1, \alpha_2 + \alpha_2^{-1} \beta_2^2 h^2 \big)}} + h \beta_1 J.$$ \qed

\section{A two-dimensional periodic medium\label{per}}

In this section we consider a sequence $\Sigma_n$ of matrix valued functions (not necessarily symmetric) in $L^\infty(\R^2)^{2 \times 2}$, which satisfies the following assumptions:

\begin{itemize}

\item[1]. $\Sigma_n$ is $Y$-periodic, where $Y:=(0,1)^2$, i.e.,\begin{eqnarray}\forall\,n \in \mathbb N, \forall\,\kappa \in \mathbb Z^2, \ \ \Sigma_n(. + \kappa) = \Sigma_n(.) \ \ \textrm{a.e.} \ \textrm{in} \ \R^2,\label{periodic1}\end{eqnarray}

\item[2]. $\Sigma_n$ is equi-coercive in $\R^2$, i.e.,\begin{eqnarray}\exists\,\alpha > 0 \ \ \ \textrm{such} \ \textrm{that} \ \ \ \forall\,n \in \mathbb N, \forall\,\xi \in \R^2, \ \ \Sigma_n \xi \cdot \xi \geq \alpha |\xi|^2 \ \ \textrm{a.e.} \ \textrm{in} \ \R^2.\label{periodic2}\end{eqnarray}
\end{itemize}

Let $\varepsilon_n$ be a sequence of positive numbers which tends to $0$. From the sequences $\Sigma_n$ and $\varepsilon_n$ we define the highly oscillating sequence of matrix-valued functions $\sigma_n$ by\begin{eqnarray}\sigma_n(x) = \displaystyle \Sigma_n\left(\frac{x}{\varepsilon_n} \right), \ \ \textrm{a.e.} \ x \in \R^2.\label{periodic3}\end{eqnarray}

\noindent By virtue of (\ref{periodic1}) and (\ref{periodic2}), $\sigma_n$ is an equi-coercive sequence of $\varepsilon_n$-periodic matrix-valued functions in $L^\infty(\R^2)^{2 \times 2}$. For a fixed $n \in \mathbb N$, let $(\sigma_n)_*$ be the constant matrix defined by\begin{eqnarray}
\forall\,\lambda, \mu \in \R^2, \quad (\sigma_n)_* \lambda \cdot \mu =\int_Y \Sigma_n \nabla W_n^\lambda \cdot \nabla W_n^\mu \ \text{d}y, \label{defsigman*}
\end{eqnarray}where, for any $\lambda \in \R^2$, $W_n^\lambda \in H^1_\sharp(Y)$, the set of $Y$-periodic functions belonging to $H^1_{loc}(\R^2)$, is the solution of the auxiliary problem

\begin{equation}
\displaystyle \int_Y \big(W_n^\lambda - \lambda \cdot y \big) \ \text{d}y =0 \quad \text{and} \quad \text{div}\big( \Sigma_n \nabla W_n^\lambda\big) = 0 \quad \text{in} \ \mathscr{D}'(\R^2)
\end{equation} or equivalently\begin{equation}
\left\{\!\! \begin{array}{r c l l}
\displaystyle &\displaystyle\int_Y \Sigma_n \nabla W_n^\lambda \cdot \nabla \varphi \ \text{d}y = 0 , \quad \forall\,\varphi \in H^1_{\sharp}(Y) \vspace{0.2cm}\\
&\displaystyle\int_Y \big(W_n^\lambda(y) - \lambda \cdot y \big) \ \mathrm{d}y = 0.
\end{array} \right. \;\;\;\label{pb4}\end{equation} Set \begin{equation}w_n^\lambda(x) := \varepsilon_n W_n^\lambda\left( \frac{x}{\varepsilon_n}\right),\quad \text{for} \ x \in \Omega,\label{w_n}\end{equation}

\noindent and \begin{equation}
w_n : = (w_n^{e_1},w_n^{e_2}) = (w_n^{1},w_n^{2}) \label{corrector1}.
\end{equation}

\subsection{A uniform convergence result}

\begin{Theo}

Let $\Omega$ be a bounded open subset of $\R^2$ with a Lipschitz boundary. Consider a highly oscillating sequence of matrix-valued functions $\sigma_n$ satisfying (\ref{periodic1}), (\ref{periodic2}), (\ref{periodic3}) and the constant matrix $(\sigma_n)_*$ defined by (\ref{defsigman*}). We assume that\begin{equation} (\sigma_n)_* \longrightarrow \sigma_* \ \textrm{in} \ \R^{2 \times 2}. \label{a1}\end{equation}

\noindent Consider, for $f \in H^{-1}(\Omega) \cap W^{-1,q}(\Omega)$ with $q>2$, the solution $u_n$ of the problem

\begin{equation}
\mathcal{P}_n \; \left\{\!\! \begin{array}{r c l l}
- \mathrm{div} (\sigma_n \nabla u_n ) & = & f & \ \ \mathrm{in} \ \Omega\\
u_n & = & 0 & \ \ \mathrm{on} \ \partial \Omega .
\end{array} \right. \;\;\;\label{pb1}\end{equation}

\noindent Then, $u_n$ converges uniformly to the solution $u \in H_0^1(\Omega)$ of

\begin{equation}
\mathcal{P} \; \left\{\!\! \begin{array}{r c l l}
- \mathrm{div} (\sigma_* \nabla u ) & = & f & \ \ \mathrm{in} \ \Omega\\
u & = & 0 & \ \ \mathrm{on} \ \partial \Omega .
\end{array} \right. \;\;\;\label{pb2}\end{equation}

\noindent Moreover we have the corrector result, with the $\varepsilon_n Y$-periodic sequence $w_n$ defined in \eqref{corrector1}: \begin{equation}
\nabla u_n - \Sum_{i=1}^2 \partial_i u \ \nabla w_n^i \longrightarrow 0 \quad \text{in} \ L^1(\Omega)^2. \label{ttt}
\end{equation}

\label{thunifconv}

\end{Theo}

\begin{Rem}

The first point of Theorem \ref{thunifconv} is an extension to the non-symmetric case of the results of \cite{bcconvunif} and \cite{bcunifconv}. The uniform convergence of $u_n$ is a straightforward consequence of Theorem 2.7 of \cite{bcunifconv} taking into account that in the present case $\sigma_n \in L^\infty(\Omega)^{2 \times 2}$ for a fixed $n$. The fact that $f \in W^{-1,q}(\Omega)$ with $q >2$ ensures the uniform convergence.
\label{remunif}
\end{Rem}

\noindent \textbf{Proof of Theorem \ref{thunifconv}}.
\par\smallskip\noindent
\textit{Derivation of the limit problem $\mathcal{P}$.}
\smallskip

We only have to show that $u$ is the solution of $\mathcal{P}$ in \eqref{pb2}. We consider a corrector \linebreak[4] \mbox{$D \widetilde{w}_n : \R^2 \longrightarrow \R^{2 \times 2}$} associated with $\sigma_n^{\text{{T}}}$ defined by\[\widetilde{w}_n(x) := \varepsilon_n \widetilde{W}_n\bigg( \frac{x}{\varepsilon_n}\bigg) = \left(\varepsilon_n \widetilde{W}^1_n\bigg( \frac{x}{\varepsilon_n}\bigg),\varepsilon_n \widetilde{W}^2_n\bigg( \frac{x}{\varepsilon_n}\bigg) \right)\]where for $i = 1,2$, $\widetilde{W}_n^i \in H^1_\sharp (Y)$ is the solution of the auxiliary problem \begin{equation}
\displaystyle \int_Y \big( \widetilde{W}^i_n - e_i \cdot x \big) \dx =0 \quad \text{and} \quad \text{div}\left( \Sigma_n^\text{T} \nabla \widetilde{W}^i_n \right) = 0 \quad \text{in} \ \mathscr D'(\R^2). \label{corrunif}
\end{equation}

Again, thanks to Theorem 2.7 of \cite{bcunifconv}, $\widetilde{w}_n$ converges uniformly to the identity in $\Omega$ by the integral condition \eqref{corrunif}. Let $\varphi \in \mathscr{D}(\Omega)$. We have, using the Einstein convention, by integrating by parts and by the Schwarz theorem ($\partial^2_{i,j} \varphi = \partial^2_{j,i} \varphi$) \begin{align*} & \hspace{0.5cm} \int_\Omega \sigma_n \nabla u_n \cdot \nabla \left(\varphi (\widetilde{w}_n) \right)\dx\\ &= \int_\Omega \nabla u_n \cdot \sigma_n^\mathrm{T} \nabla \widetilde{w}_n^i (\partial_i \varphi) (\widetilde{w}_n) \dx \\
&= \underbrace{\int_\Omega \sigma_n^\mathrm{T} \nabla \widetilde{w}_n^i \cdot \nabla (u_n\partial_i \varphi (\widetilde{w}_n)) \dx}_{= 0 } - \int_\Omega \sigma_n^\mathrm{T} \nabla \widetilde{w}_n^i \cdot \nabla \widetilde{w}_n^j \ \ \partial^2_{i,j} \varphi (\widetilde{w}_n) \ u_n\dx \\
&= - \int_\Omega \sigma_n \nabla \widetilde{w}_n^i \cdot \nabla \widetilde{w}_n^i \ \ \partial^2_{i,i} \varphi (\widetilde{w}_n) \ u_n \dx - \int_\Omega \sigma_n^\text{T} \nabla \widetilde{w}_n^2 \cdot \nabla \widetilde{w}_n^1 \ \ \partial^2_{2,1} \varphi (\widetilde{w}_n) \ u_n \dx\\
& \textcolor{white}{espace} - \int_\Omega \sigma_n^\mathrm{T}\nabla \widetilde{w}_n^1 \cdot \nabla \widetilde{w}_n^2 \ \ \partial^2_{1,2} \varphi (\widetilde{w}_n) \ u_n \dx \\
&= - \int_\Omega \sigma_n \nabla \widetilde{w}_n^i \cdot \nabla \widetilde{w}_n^i \ \ \partial^2_{i,i} \varphi (\widetilde{w}_n) \ u_n \dx - \int_\Omega \sigma_n \nabla \widetilde{w}_n^1 \cdot \nabla \widetilde{w}_n^2 \ \ \partial^2_{1,2} \varphi (\widetilde{w}_n) \ u_n \dx\\
& \textcolor{white}{espace} - \int_\Omega \sigma_n^\mathrm{T}\nabla \widetilde{w}_n^1 \cdot \nabla \widetilde{w}_n^2 \ \ \partial^2_{1,2} \varphi (\widetilde{w}_n) \ u_n \dx \\
&= - \int_\Omega \sigma_n^s \nabla \widetilde{w}_n^i \cdot \nabla \widetilde{w}_n^i \ \ \partial^2_{i,i} \varphi (\widetilde{w}_n) \ u_n \dx - 2 \int_\Omega \sigma_n^s \nabla \widetilde{w}_n^1 \cdot \nabla \widetilde{w}_n^2 \ \ \partial^2_{1,2} \varphi (\widetilde{w}_n) \ u_n \dx. \end{align*}This leads us to the equality\begin{equation} \langle f, \varphi (\widetilde{w}_n)\rangle_{H^{-1}(\Omega),H_0^1(\Omega)} = \int_\Omega \sigma_n \nabla u_n \cdot \nabla \left(\varphi (\widetilde{w}_n) \right)\dx = - \int_\Omega \sigma_n^s \nabla \widetilde{w}_n^i \cdot \nabla \widetilde{w}_n^j \ \ \partial^2_{i,j} \varphi (\widetilde{w}_n) \ u_n\dx. \label{IPP}\end{equation}

To study the convergence of the last term of \eqref{IPP}, we first show that $\sigma_n^s \nabla \widetilde{w}_n^i \cdot \nabla \widetilde{w}_n^j$ is bounded in $L^1(\Omega)$. We have, by periodicity and the Cauchy-Schwarz inequality \begin{align*} \displaystyle \int_\Omega |\sigma_n^s \nabla \widetilde{w}_n^i \cdot \nabla \widetilde{w}_n^j| \dx &= \displaystyle \int_\Omega |\Sigma_n^s \nabla \widetilde{W}_n^i \cdot \nabla \widetilde{W}_n^j| \left(\frac{x}{\varepsilon_n} \right) \dx \\
&\leq C \int_{Y} |\Sigma_n^s \nabla \widetilde{W}_n^i \cdot \nabla \widetilde{W}_n^j| \dx \\
&\leq C \sqrt{ \int_{Y} \big|\Sigma_n^s \nabla \widetilde{W}_n^i \cdot \nabla \widetilde{W}_n^i\big| \dx} \sqrt{ \int_{Y} \big|\Sigma_n^s \nabla \widetilde{W}_n^j \cdot \nabla \widetilde{W}_n^j\big| \dx} \\
&\leq C \sqrt{(\sigma_n)_* e_i \cdot e_i} \sqrt{(\sigma_n)_* e_j \cdot e_j}
\end{align*}which is bounded by the hypothesis \eqref{a1}. Therefore,\begin{equation}\sigma_n^s \nabla \widetilde{w}_n^i \cdot \nabla \widetilde{w}_n^j \ \text{is bounded in} \ L^1(\Omega). \label{boundsympart} \end{equation}

Due to the periodicity, we know that for $i,j = 1,2$, \begin{align*} 2 \sigma_n^s \nabla \widetilde{w}_n^i \cdot \nabla \widetilde{w}_n^j = \sigma_n^\text{T} \nabla \widetilde{w}_n^i \cdot \nabla \widetilde{w}_n^j + \sigma_n^\text{T} \nabla \widetilde{w}_n^j \cdot \nabla \widetilde{w}_n^i \rightharpoonup \left(\sigma_*\right)^\text{T} e_i \cdot e_j + \left(\sigma_*\right)^\text{T} e_j \cdot e_i = 2 \left(\sigma_*\right)^s e_i \cdot e_j \end{align*}weakly-$*$ in $\mathcal M(\Omega)$. Hence, we get that\begin{equation} \sigma_n^s \nabla \widetilde{w}_n^i \cdot \nabla \widetilde{w}_n^j \rightharpoonup \left(\sigma_*\right)^s e_i \cdot e_j \quad \text{weakly-$*$ in} \ \mathcal M (\Omega).\label{convpsym} \end{equation}Moreover, $ \partial^2_{i,j} \varphi(\widetilde{w}_n) \ u_n $ converges uniformly to $\partial^2_{i,j} \varphi \ u$. Thus, by passing to the limit in \eqref{IPP}, we have, again with the Einstein convention\[ \langle f, \varphi \rangle_{H^{-1}(\Omega),H_0^1(\Omega)} = - \int_{\Omega} \left(\sigma_*\right)^s e_i \cdot e_j \ \ \partial_{i,j}^2 \varphi \ u \dx = - \int_\Omega \sigma_* : \nabla^2 \varphi \ u \dx.\]

\noindent Therefore, by integrating by parts and using $\varphi = 0$ on $\partial \Omega$, \begin{equation} \int_\Omega \sigma_* \nabla u \cdot \nabla \varphi \dx = \langle f, \varphi \rangle_{H^{-1}(\Omega),H_0^1(\Omega)}. \label{IPPfin2}\end{equation}

\bigskip

\noindent \textit{Proof of the corrector result}

\smallskip

First of all, we show that the corrector function $w_n$ is bounded in $H^1(\Omega)^2$. By the definition \eqref{corrector1} of $w_n$, the $Y$-periodicity of $W_n^{e_i}$ and the equi-coercivity of $\Sigma_n$, we have, for $i=1,2$,\begin{equation} \alpha \ ||\nabla w_n^i||_{L^2(\Omega)^2}^2 \leq C \alpha \ ||\nabla W_n^{e_i}||_{L^2(Y)^2}^2 \leq C \int_Y \Sigma_n \nabla W_n^i \cdot \nabla W_n^i \dx = C \ (\sigma_n)_* e_i \cdot e_i \label{boundwn}\end{equation} which is bounded. This inequality combined with the uniform convergence of $w_n$ yields to the boundedness of $w_n$ in $H^1(\Omega)^2$.

Let us consider an approximation $u^\delta \in \mathscr D(\Omega)$ of $u$ such that\begin{equation}
|| u - u^\delta||_{H_0^1(\Omega)} \leq \delta.
\label{approxgrad}
\end{equation}

On the one hand, we have \begin{align*}\int_\Omega \sigma_n \nabla u_n \cdot \nabla \big(u_n - u^\delta(w_n) \big)\dx &= \langle f, \big(u_n - u^\delta(w_n) \big) \rangle_{H^{-1}(\Omega),H_0^1(\Omega)} . \end{align*} Since $w_n$ converges uniformly to identity on $\Omega$ and is bounded in $H^1(\Omega)$ (see \eqref{boundwn}), with $u^\delta \in \mathscr D(\Omega)$, $u^\delta(w_n)$ converges weakly to $u^\delta$ in $H_0^1(\Omega)$. Hence, by the weak convergence of $u_n$ to $u$ in $H_0^1(\Omega)$ and \eqref{approxgrad}, we can pass to the limit the previous inequality and obtain, for any $\delta > 0$,\begin{equation}\limsup\limits_{n \to \infty}\left|\int_\Omega \sigma_n \nabla u_n \cdot \nabla \big(u_n - u^\delta(w_n) \big) \dx \right| = \left| \langle f, u - u^\delta \rangle_{H^{-1}(\Omega),H_0^1(\Omega)} \right| \leq C \delta \label{tend1}.\end{equation}

On the other hand, similarly to the proof of the first point \eqref{IPP}, we are led to the equality\begin{equation} \int_\Omega \sigma_n \nabla \big(u^\delta(w_n)\big) \cdot \nabla \big(u_n - u^\delta(w_n) \big) \dx = - \int_\Omega \sigma_n^s \nabla w_n^i \cdot \nabla w_n^j \ \ \partial^2_{i,j} u^\delta(w_n) \ \big(u_n - u^\delta(w_n) \big) \dx. \label{IPP2}\end{equation}

As in the first point, $\sigma_n^s \nabla w_n^i \cdot \nabla w_n^j$ is bounded in $L^1(\Omega)$ (see \eqref{boundsympart}), $u_n$ converges uniformly to $u$ and $\partial_{i,j} u^\delta(w_n)$ converges uniformly to $\partial_{i,j} u^\delta$ because $u^\delta$ is a $\mathscr D(\Omega)$ function. By passing to the limit in \eqref{IPP2} \begin{equation} \int_\Omega \sigma_n \nabla \big(u^\delta(w_n)\big) \cdot \nabla \big(u_n - u^\delta(w_n) \big) \dx \underset{n \to \infty}{\longrightarrow} - \int_\Omega ( \sigma_*)^s e_i \cdot e_j \ \ \partial^2_{i,j} u^\delta \ \big(u - u^\delta \big) \dx. \label{tend2} \end{equation}

\noindent Moreover, like in \eqref{IPPfin2} we have \begin{equation}
\int_\Omega ( \sigma_*)^s e_i \cdot e_j \ \ \partial^2_{i,j} u^\delta \ \big(u - u^\delta \big) \dx = \int_\Omega \sigma_* \nabla u^\delta \cdot \nabla \big(u - u^\delta \big) \dx.
\end{equation}

\noindent By combining this equality with the convergence \eqref{tend2}, we obtain the inequality \begin{align}\lim\limits_{n \to \infty}\left|\int_\Omega \sigma_n \nabla \big(u^\delta(w_n)\big) \cdot \nabla \big(u_n - u^\delta(w_n) \big) \dx \right| &\leq \left| \int_\Omega \sigma_* \nabla u^\delta \cdot \nabla \big(u - u^\delta \big)\right| \\
&\leq C |\sigma_*| \ || \nabla u^\delta||_{L^2(\Omega)^2} || \nabla \big( u - u^\delta \big)||_{L^2(\Omega)^2}\leq C \delta \label{tend4}.\end{align}

\noindent Thus, by adding \eqref{tend1} and \eqref{tend4}, we have \begin{equation*}
\limsup\limits_{n \to \infty} \int_\Omega \sigma_n \nabla \big(u_n - u^\delta(w_n) \big) \cdot \nabla \big(u_n - u^\delta(w_n) \big) \dx \leq C \delta
\label{majomax}
\end{equation*} which leads us, by equi-coercivity, to\begin{equation} \begin{array}{r l} &\displaystyle \hspace{0.4cm} \limsup \limits_{n \to \infty} \ \alpha \ ||\nabla(u_n - u^\delta(w_n))||^2_{L^2(\Omega)^2} \vspace{0.2cm}\\
&  \leq \displaystyle \limsup \limits_{n \to \infty} \left|\int_\Omega \sigma_n \nabla \big(u_n - u^\delta(w_n) \big) \cdot \nabla \big(u_n - u^\delta(w_n) \big) \dx \right| \leq C \delta. \label{final}\end{array}\end{equation}

Thus, the Cauchy-Schwarz inequality, the boundedness of $\nabla w_n^i$ in $L^2(\Omega)^2$ \eqref{boundwn} and the Einstein convention give, for any $\delta >0$, \begin{align*} & \hspace{0.5cm} ||\nabla u_n - \nabla w_n^i \partial_i u||_{L^1(\Omega)^2}\\
 & \leq ||\nabla u_n - \nabla w_n^i \partial_i u^\delta ||_{L^1(\Omega)^2} + ||\nabla w_n^i \partial_i \big(u^\delta - u\big)||_{L^1(\Omega)^2} \\
& \leq ||\nabla u_n - \nabla w_n^i \partial_i u^\delta ||_{L^1(\Omega)^2} + ||\nabla w_n^i||_{L^2(\Omega)^2} ||\partial_i \big(u^\delta - u\big)||_{L^2(\Omega)} \\
&\leq ||\nabla u_n - \nabla w_n^i \partial_i u^\delta ||_{L^1(\Omega)^2} + C \delta \\
& \leq ||\nabla u_n - \nabla w_n^i \partial_i u^\delta(w_n) ||_{L^1(\Omega)^2} +||\nabla w_n^i \big(\partial_i u^\delta - \partial_i u^\delta(w_n)\big)||_{L^1(\Omega)^2} + C \delta \\
& \leq ||\nabla u_n - \nabla w_n^i \partial_i u^\delta(w_n) ||_{L^1(\Omega)^2} +||\nabla w_n^i|| _{L^2(\Omega)^2} || \partial_i u^\delta - \partial_i u^\delta(w_n)||_{L^2(\Omega)} + C \delta \\
& \leq ||\nabla u_n - \nabla w_n^i \partial_i u^\delta(w_n) ||_{L^1(\Omega)^2} + C || \partial_i u^\delta - \partial_i u^\delta(w_n)||_{L^2(\Omega)} + C \delta.
\end{align*}Since $u^\delta \in \mathscr D(\Omega)$ and $w_n$ converges uniformly to the identity on $\Omega$, the second term of the last inequality converges to 0. Hence, we get that \begin{equation}\limsup \limits_{n \to \infty} ||\nabla u_n - \nabla w_n^i \partial_i u||_{L^1(\Omega)^2} \leq \limsup \limits_{n \to \infty} ||\nabla u_n - \nabla w_n^i \partial_i u^\delta(w_n) ||_{L^1(\Omega)^2} + C \delta.
\end{equation} Finally, this inequality combined with \eqref{final} gives, for any $\delta > 0$, \begin{equation*}
0 \leq \limsup \limits_{n \to \infty} ||\nabla u_n - \nabla w_n^i \partial_i u||_{L^1(\Omega)^2} \leq C \sqrt \delta + C \delta,
\end{equation*} which implies the corrector result \eqref{ttt}. \qed

\begin{Rem}

If the solution $u$ is a $\mathscr C^2$ function, then the convergence \eqref{ttt} holds true in $L_{\text{loc}}^2(\Omega)$ since we may take $u = u^\delta$.

\end{Rem}

\subsection{A two-phase result \label{2phasres}}

Here, we recall a two-phase result due to G.W. Milton (see \cite{MilBook} pp. 61--65) using the Dykhne transformation.

In order to apply the previous theorem, we reformulate Milton's calculus in such a way that every coefficient depends on $n$. We then consider, for a fixed $n$, the periodic homogenization of a conductivity $\sigma_n(h)$ to obtain $(\sigma_n)_*(h)$ through the link between the homogenization of the transformed conductivity and $(\sigma_n)_*(h)$ given by formula (4.16) in \cite{MilBook}. Finally, we study the limit  of $(\sigma_n)_*(h)$ through the asymptotic behavior of the coefficients of the transformation, and apply Theorem \ref{thunifconv} in the example Section \ref{example}.

In this section we consider a two-phase periodic isotropic medium. Let $\chi_n$ be a sequence of characteristic functions of subsets of $Y$. We define for any $\alpha_1 >0$, $\beta_1 \in \R$, any sequences $\alpha_{2,n} >0$, $\beta_{2,n} \in \R$ and any $h \in \R$, a parametrized conductivity $\Sigma_n(h)$: \begin{equation}\Sigma_n(h) = (1-\chi_n) (\alpha_1 I_2 + h \beta_1 J) + \chi_n (\alpha_{2,n} I_2 + h \beta_{2,n}J) \quad \text{in} \ Y. \label{2phas1}\end{equation} We still denote by $\Sigma_n(h)$ the periodic extension to $\R^2$ of $\Sigma_n(h)$ (which satisfies \eqref{periodic1}). We assume that $\Sigma_n(h)$ satisfies \eqref{periodic2}, and define $\sigma_n(h)$ by \eqref{periodic3} and $(\sigma_n)_*(h)$ by \eqref{defsigman*}.

\bigskip

We have the following result based on an analysis of \cite{MilBook} (pp. 61--65).

\begin{Prop}

Let $\chi_n$ be a sequence of characteristic functions of subsets of $Y$, $\alpha_1, \alpha_2 > 0$,  a positive sequence $\alpha_{2,n}$, $\beta_1, \beta_2 \in \R$, and a sequence $\beta_{2,n}$  such that\begin{equation}
\lim \limits_{n \to \infty} \alpha_{2,n} = \infty , \quad \liminf \limits_{n \to \infty} |\beta_{2,n} - \beta_1|  > 0, \quad \text{and} \quad \lim \limits_{n \to \infty}\frac{\beta_{2,n}}{\alpha_{2,n}} = \frac{\beta_2}{\alpha_2}.\label{g1}\end{equation}Assume that the effective conductivity in the absence of a magnetic field\begin{equation}
\left(\sigma_n^0\right)_*(\gamma_{1,n},\gamma_{2,n}) \label{sn} \quad \text{is bounded when}
\quad\lim_{n\to\infty}\gamma_{1,n}=\alpha_1\;\;\mbox{and}\;\;\lim_{n\to\infty}{\frac{\gamma_{2,n}}{\alpha_{2,n}}}=\gamma_2>0.
\end{equation}

\noindent Then, there exist two parametrized positive sequences $\alpha'_{1,n}(h), \alpha'_{2,n}(h)$ such that\begin{equation}\lim \limits_{n \to \infty} \alpha'_{1,n}(h) = \alpha_1 \quad \textrm{and} \quad \alpha'_{2,n}(h) \underset{n \to \infty}{\sim} \frac{\alpha_2^2 + h^2 \beta_2^2}{\alpha_2^2}\,\alpha_{2,n},\label{i44} \end{equation}and
\begin{equation} (\sigma_n)_*(h) = \left(\sigma_n^0\right)_*\big(\alpha'_{1,n}(h),\alpha'_{2,n}(h)\big) + h \beta_1 J + \underset{n \to \infty}{o}(1)\end{equation}where $\left(\sigma^0_n\right)_*\big(\alpha'_{1,n}(h),\alpha'_{2,n}(h)\big)$ is bounded.
\label{Mil}
\end{Prop}

\begin{Rem}
In view of condition \eqref{g1}, the case where $\beta_{2,n}$ tends to $\beta_1$ corresponds to perturb the symmetric conductivity \begin{equation*} \sigma_n^s = (1 - \chi_n) \alpha_1 I_2 + \chi_n \alpha_{2,n} I_2\end{equation*} by \begin{equation*}\sigma_n^s + \beta_1 J + \underset{n \to \infty }{o}(1).\end{equation*} Then it is clear that \begin{equation*} (\sigma_n)_*(h) = \big(\sigma_n^s\big)_* + \beta_1 J + \underset{n \to \infty }{o}(1). \end{equation*}
\end{Rem}

\noindent \textbf{Proof of Proposition \ref{Mil}}. The proof is divided into two parts. After applying Milton's computation (pp. 61--64 of \cite{MilBook}), we study the asymptotic behavior of the different coefficients.

\textcolor{white}{t} 

\noindent \textit{First step}: Applying Dykhne's transformation through Milton's computations.

In order to make the Dykhne's transformation following Milton \cite{MilBook} (pp. 62--64), we consider two real coefficients $a_n$ and $b_n$ such that\begin{equation}
\sigma'_n := \big( a_n \sigma_n(h) + b_n J\big)\big(a_n I_2 + J \sigma_n(h) \big)^{-1} = a_n \big( \sigma_n(h) + (a_n)^{-1} b_n J\big)\big( a_n I_2 + J \sigma_n(h) \big)^{-1}
\label{k1}\end{equation} is symmetric and, more precisely, according to Notation \ref{defabs}, reads as \begin{equation} \sigma'_n = (1 - \chi_n)\alpha'_{1,n}(h) I_2 + \chi_n \alpha'_{2,n}(h) I_2 = \sigma_n^0\big(\alpha'_{1,n}(h),\alpha'_{2,n}(h)\big) \label{zero}. \end{equation}

Then, using the complex representation\begin{equation}
\alpha I_2 + \beta J \longleftrightarrow \alpha + \beta i
\end{equation} suggested by Tartar \cite{comTar}, the constants $a_n$, $b_n$ must satisfy\begin{equation} \alpha'_{1,n}(h) = \frac{a_n(\alpha_{1} + i h \beta_1) + i b_n}{a_n + i(\alpha_{1} + i h \beta_1)} \in \R \quad \textrm{and} \quad \alpha'_{2,n}(h) = \frac{a_n(\alpha_{2,n} + i h \beta_{2,n}) + i b_n}{a_n + i(\alpha_{2,n} + i h \beta_{2,n})} \in \R, \end{equation} which implies that \begin{equation}
b_n = \frac{-a_n^2 h \beta_1 + a_n \Delta_1}{a_n - h \beta_1} = \frac{-a_n^2 h \beta_{2,n} + a_n \Delta_{2,n}}{a_n - h \beta_{2,n}} \label{g9}.
\end{equation}Denoting $\Delta_1 := \alpha_1^2 + h^2 \beta_1^2$ and $\Delta_{2,n} := \alpha_{2,n}^2 + h^2 \beta_{2,n}^2$ (thanks to \eqref{g1}, $n$ is considered to be larger enough such that $\beta_{2,n} - \beta_1 \neq 0$ and $a_n$ is real), the equality \eqref{g9} provides two non-zero solutions for $a_n$:\begin{equation}
\displaystyle a_n = \frac{\Delta_{2,n} - \Delta_1 +  \sqrt{(\Delta_{2,n} - \Delta_1)^2 + 4h^2(\beta_{2,n} - \beta_1)(\beta_{2,n} \Delta_1 - \beta_1 \Delta_{2,n})}}{2h(\beta_{2,n} - \beta_1)}, \label{g2}\end{equation}and \begin{equation*}
\displaystyle a_n^- = \frac{\Delta_{2,n} - \Delta_1 -  \sqrt{(\Delta_{2,n} - \Delta_1)^2 + 4h^2(\beta_{2,n} - \beta_1)(\beta_{2,n} \Delta_1 - \beta_1 \Delta_{2,n})}}{2h(\beta_{2,n} - \beta_1)}. \end{equation*}The value \eqref{g2} is associated with a positive matrix $\sigma'_n$, while $a_n^-$ leads us to the negative matrix $\sigma^-_n = - J( \sigma'_n)^{-1} J^{-1}$ to exclude (see \cite{MilHall} for more details).

\textcolor{white}{t}

\noindent \textit{Second step}: asymptotic behavior of the coefficients and the homogenized matrix.

One the one hand, by the equality \eqref{g2} combined with \eqref{g1}, we have\begin{equation*}
\lim \limits_{n \to \infty} a_n \frac{h (\beta_{2,n} - \beta_1)}{ \alpha^2_{2,n}} = \frac{\alpha_2^2 + h^2 \beta_2^2}{\alpha_2^2}
\end{equation*}
which clearly implies that\begin{equation}
 a_n \underset{n \to \infty}{\sim} \frac{\alpha_2^2 + h^2 \beta_2^2}{ \alpha_2^2} \frac{\alpha_{2,n}^2}{h(\beta_{2,n} - \beta_1)} \quad \text{and} \quad a_n - h \beta_{2,n} \underset{n \to \infty}{\sim}\frac{\alpha_{2,n}^2}{h(\beta_{2,n} - \beta_1)}. \label{g3}
\end{equation}

On the other hand, \eqref{g1}, \eqref{g3} and the first equality of \eqref{g9} give\begin{equation}
b_n = -a_n h \beta_1 + \Delta_1 + \underset{n \to \infty}{o}(1). \label{g22}
\end{equation}

From \eqref{g1}, \eqref{g2}, \eqref{g3} and \eqref{g22} we deduce the following asymptotic behavior for the modified phases:\begin{equation}\lim \limits_{n \to \infty} \alpha'_{1,n}(h) = \alpha_1 \quad \text{and} \quad \lim \limits_{n \to \infty} \frac{\alpha'_{2,n}(h)}{\alpha_{2,n}} = \frac{\alpha_2^2 + h^2 \beta_2^2}{\alpha_2^2}. \label{2phas}\end{equation}

To consider $\big(\sigma'_n \big)_*$, we need to verify that $\sigma'_n$ is equi-coercive. We have, by denoting for any $\xi \in \R^2$, $\nu_n = \big(a_n I_2 + J \sigma_n(h)\big)^{-1} \xi$, \begin{equation*}
\forall\,\ \xi \in \R^2, \quad \sigma'_n \xi \cdot \xi = \big( a_n \sigma_n(h) + b_n J\big) \nu_n \cdot \big(a_n I_2 + J \sigma_n(h)\big) \nu_n = (a_n^2 + b_n ) \sigma_n(h) \nu_n \cdot \nu_n
\end{equation*} and, because $a_n^{-1} \sigma_n(h)$ is bounded in $L^\infty(\Omega)^{2 \times 2}$ by \eqref{g3}, \begin{equation*}
\forall\,\ \xi \in \R^2, \quad |\xi| = \big| a_n \nu_n + J \sigma_n(h) \nu_n \big| \leq a_n (1+C) |\nu_n| . 
\end{equation*} The equi-coercivity of $\sigma_n(h)$ gives \begin{equation}
\exists\,C > 0, \quad \forall\,\xi \in \R^2, \quad \sigma'_n \xi \cdot \xi \geq \frac{C}{(1+C)^2} \frac{a_n^2 + b_n}{a_n^2} |\xi|^2 
\end{equation}
that is, for $n$ larger enough, by \eqref{g3} and \eqref{g22}, $\sigma'_n$ is equi-coercive.

We can now apply the Keller-Dykhne duality theorem (see, e.g., \cite{Keller,Dyk}) to equality \eqref{k1} to obtain\begin{equation}
(\sigma'_n)_* = \big( a_n (\sigma_n)_* + b_n J\big)\big( a_n I_2 + J (\sigma_n)_* \big)^{-1}. \label{sigman'*}
\end{equation} Moreover, by inverting this transformation, we have \[(\sigma_n)_*(h) = \big( a_n I_2 - (\sigma'_n)_* J  \big)^{-1}\big( a_n (\sigma'_n)_* - b_n J\big).\] Considering \eqref{g1}, \eqref{g3}, \eqref{g22}, and the boundedness of $(\sigma'_n)_*$ (as a consequence of the bound \eqref{sn}) we get that
\begin{equation} (\sigma_n)_*(h) = (\sigma'_n)_* - \frac{b_n}{a_n} J + \underset{n \to \infty}{o}(1) = (\sigma'_n)_* + h \beta_1 J + \underset{n \to \infty}{o}(1), \end{equation}which concludes the proof taking into account \eqref{zero}. \qed

\bigskip

To derive the limit of $\left(\sigma^0_n\right)_* (\alpha_{1,n}'(h), \alpha_{2,n}'(h))$, we need more information on the geometry of the high conductive phase.
To this end, we study the following example.

\subsection{A cross-like thin structure\label{example}}

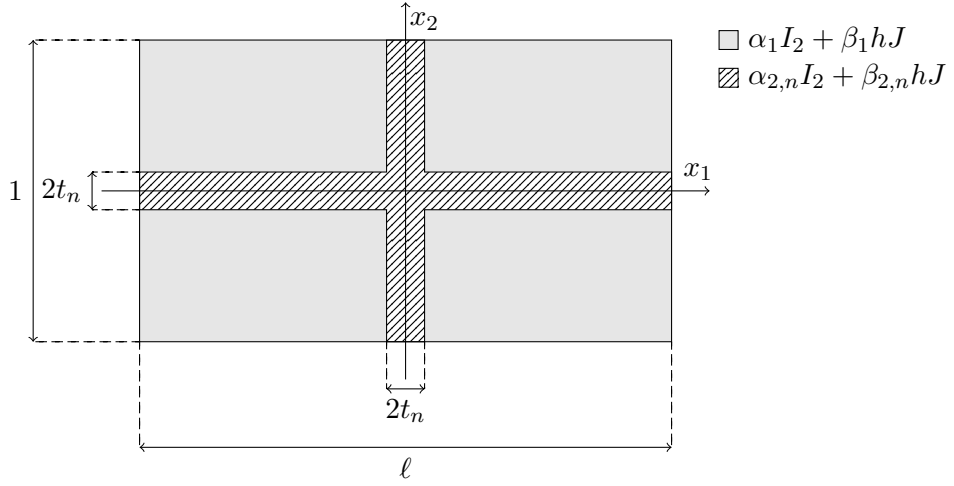
\begin{figure}[!t]
\centering
\begin{tikzpicture}[scale=0.5,xmin=-18,xmax=13.5,ymin=-6,ymax=6]
\draw (-14.2,0) node[right, color = white]{$a$};
\fill [color=gray!20 ]
(8.25,3.75) -- (8.25,4.25) -- (8.75,4.25) -- (8.75,3.75) -- cycle;
\draw (8.25,3.75) -- (8.25,4.25) -- (8.75,4.25) -- (8.75,3.75) -- cycle;
\fill [ pattern= north east lines]
(8.25,2.75) -- (8.25,3.25) -- (8.75,3.25) -- (8.75,2.75) -- cycle;
\draw (8.25,2.75) -- (8.25,3.25) -- (8.75,3.25) -- (8.75,2.75) -- cycle;
\draw (8.75,4) node[right]{$\alpha_1 I_2 + \beta_1 h J$};
\draw (8.75,3) node[right]{$\alpha_{2,n} I_2 + \beta_{2,n} h J$};
\fill [color=gray!20 ]
(-7,-4) -- (-7,-0.5) -- (-0.5,-0.5) -- (-0.5,-4) -- cycle;
\fill [color=gray!20 ]
(7,4) -- (7,0.5) -- (0.5,0.5) -- (0.5,4) -- cycle;
\fill [color=gray!20 ]
(7,-4) -- (7,-0.5) -- (0.5,-0.5) -- (0.5,-4) -- cycle;
\fill [color=gray!20 ]
(-7,4) -- (-7,0.5) -- (-0.5,0.5) -- (-0.5,4) -- cycle;
\fill [color= blue!20, pattern= north east lines]
(-7,0.5) -- (-0.5,0.5) -- (-0.5,4) -- (0.5,4) -- (0.5,0.5) -- (7,0.5) -- (7,-0.5) -- (0.5,-0.5) -- (0.5,-4) -- (-0.5,-4) -- (-0.5,-0.5) -- (-7,-0.5) -- (-7,0.5) -- cycle;
\draw (-7,0.5) -- (-0.5,0.5) -- (-0.5,4) -- (0.5,4) -- (0.5,0.5) -- (7,0.5) -- (7,-0.5) -- (0.5,-0.5) -- (0.5,-4) -- (-0.5,-4) -- (-0.5,-0.5) -- (-7,-0.5) -- (-7,0.5) -- cycle;
\draw[<->] (-0.5,-5.25) -- (0.5,-5.25)
node[midway, below]{$2 t_{n}$};
\draw[<->] (-8.25,-0.5) -- (-8.25,0.5)
node[midway, left]{$2 t_{n}$};
\draw[<->] (-7,-6.8) -- (7,-6.8)
node[midway, below]{$\ell$};
\draw[<->] (-9.8,4) -- (-9.8,-4)
node[midway, left]{$1$};
\draw[->] (-8,0) -- (8,0)
node[midway, below]{};
\draw[->] (0,-5) -- (0,5)
node[midway, below]{};
\draw (7.7,0.5) node{$x_1$};
\draw (0.5,4.5) node{$x_2$};
\draw (-7,-4) -- (-7,4) -- (7,4) -- (7,-4) -- cycle;
\draw[dash pattern = on 1mm off 1mm] (-0.5,-4) -- (-0.5,-5.25) -- cycle;
\draw[dash pattern = on 1mm off 1mm] (0.5,-4) -- (0.5,-5.25) -- cycle;
\draw[dash pattern = on 1mm off 1mm] (-7,-0.5) -- (-8.25,-0.5) -- cycle;
\draw[dash pattern = on 1mm off 1mm] (-7,0.5) -- (-8.25,0.5) -- cycle;
\draw[dashed] (-7,-4) -- (-7,-6.8) -- cycle;
\draw[dashed] (7,-4) -- (7,-6.8) -- cycle;
\draw[dashed] (-7,-4) -- (-9.8,-4) -- cycle;
\draw[dashed] (-7,4) -- (-9.8,4) -- cycle;
\end{tikzpicture}
\caption{The period of the cross-like thin structure}
\label{cross}
\end{figure}

We consider a bounded open subset $\Omega$ of $\R^2$ with a Lipschitz boundary, a real sequence $\varepsilon_n$ converging to $0$, and $f \in H^{-1}(\Omega) \cap W^{-1,q}(\Omega)$ with $q>2$. We define, for any $h \in \R$, $\alpha_1, \beta_1 >0$ and positive sequences $t_n \in (0,1/2]$, $\alpha_{2,n}$, $\beta_{2,n}$, a parametrized matrix-valued function $\Sigma_n(h)$ from the unit rectangular cell period $Y:=(-{\frac{\ell}{2}},{\frac{\ell}{2}}) \times (
-{\frac{1}{2}},{\frac{1}{2}})$, with $\ell\geq 1$, to $\R^{2 \times 2}$, by (cf. figure \ref{cross}) \begin{eqnarray}\Sigma_n(h) := \; \left\{\!\! \begin{array}{c l}
\alpha_{2,n} I_2 + \beta_{2,n} h J & \ \ \mathrm{in} \ \omega_n := \{ (x_1,x_2) \in Y \ | \ | x_1 |, |x_2| \leq t_n \} \\
\alpha_{1} I_2 + \beta_{1} h J & \ \ \mathrm{in} \ Y \setminus \omega_n
\end{array} \right. \;\;\; \label{yyyy}\end{eqnarray} Denoting again by $\Sigma_n(h)$ its periodic extension to $\R^2$, we finally consider the conductivity \begin{eqnarray}\sigma_n(h)(x) = \Sigma_n(h) \left( \frac{x}{\varepsilon_n}\right), \quad x \in \Omega,\label{defper}\end{eqnarray} and the associated homogenization problem:\begin{eqnarray}
\mathcal{P}_ {n} \; \left\{\!\! \begin{array}{r c l l}
- \mathrm{div} \big(\sigma_n(h) \nabla u_n \big) & = & f & \ \ \mathrm{in} \ \Omega\\
u_n & = & 0 & \ \ \mathrm{on} \ \partial \Omega.
\end{array} \right. \;\;\;\label{pbexplicit}
\end{eqnarray}

By virtue of Theorem \ref{thunifconv} and Proposition \ref{Mil}, we focus on the study of the limit of\linebreak[4] $\left(\sigma^0_n\right)_* \big(\alpha_{1,n}'(h), \alpha_{2,n}'(h)\big)$.

\begin{Prop}

Let $\sigma_n(h)$ be the conductivity defined by \eqref{yyyy} and \eqref{defper} and its homogenization problem \eqref{pbexplicit}. We assume that:\begin{equation}2
t_n (\ell + 1) \alpha_{2,n} \underset{n\to \infty}{\longrightarrow} \alpha_2 > 0 \ \ \textrm{and} \ \ 2
t_n (\ell + 1) \beta_{2,n} \underset{n\to \infty}{\longrightarrow} \beta_2 > 0.\label{condalphasim!}\end{equation}

\noindent Then, the homogenized conductivity is given by \[\sigma_*(h) = {\begin{pmatrix} \displaystyle{\alpha_1 + \frac{\alpha^2_2 + \beta^2_2 h^2}{( \ell + 1) \alpha_2}}  & - h \beta_1 \\ h \beta_1 & \displaystyle{\alpha_1+ \frac{\alpha^2_2 + \beta^2_2 h^2}{\ell ( \ell + 1) \alpha_2}} \end{pmatrix}}.\]
\label{thcross}

\end{Prop}

\begin{Rem}
The previous proposition does not respect exactly the framework defined at the beginning of this section because the period cell is not the unit square $Y = (0,1)^2$: we can nevertheless extend all this section to any type of period cells.
\end{Rem}

\begin{Rem}
The condition \eqref{condalphasim!} is a condition of boundedness in $L^1(\Omega)^{2 \times 2}$ of $\sigma_n$ because \[|\omega_n| = 2 t_n (\ell + 1) - 4 t_n^2 \sim 2 t_n (\ell + 1),\] which will ensure the convergence of $\left(\sigma^0_n\right)_*$.
\end{Rem}

\bigskip

\noindent \textbf{Proof of Proposition \ref{thcross}}. In order to apply Proposition \ref{Mil}, we consider two positive sequences $\alpha'_{1,n}(h), \alpha'_{2,n}(h)$ satisfying\begin{equation}\lim \limits_{n \to \infty} \alpha'_{1,n}(h) = \alpha_1 \quad \textrm{and} \quad \alpha'_{2,n}(h) \underset{n \to \infty}{\sim} \frac{\alpha_2^2 + h^2 \beta_2^2}{\alpha_2^2}\,\alpha_{2,n}.\label{55555}\end{equation}
We will study the homogenization of $\sigma'_n:=\sigma^0_n\big(\alpha'_{1,n}(h), \alpha'_{2,n}(h)\big).$
\par
To this end, consider a corrector $W_n^\lambda = \lambda \cdot x - X_n^\lambda$ in the Murat-Tartar sense (see, e.g., \cite{Mur2}) associated with \begin{eqnarray}\Sigma'_n := \; \left\{\!\! \begin{array}{c l}
\alpha'_{2,n}(h) \ I_2 & \ \ \mathrm{in} \ \omega_n = \{ (x_1,x_2) \in Y \ | \ | x_1 |, |x_2| \leq t_n \} \\
\alpha'_{1,n}(h)  \ I_2 & \ \ \mathrm{in} \ Y \setminus \omega_n
\end{array} \right. \;\;\; \end{eqnarray} and defined by\begin{equation} \; \left\{\!\! \begin{array}{l}
\vspace{0.2cm}\mathrm{div} \big(\Sigma'_n \nabla X_n^\lambda \big) = \mathrm{div} \big(\Sigma'_n \lambda \big) \ \ \mathrm{in} \ \mathscr D'(\R^2)\\
\vspace{0.2cm}X_n^\lambda \ \ \mathrm{is} \ Y-\textrm{periodic} \\
\displaystyle \int_Y X_n^\lambda \ \textrm{d}y = 0.
\end{array} \right. \;\;\;\label{defcorrsym}\end{equation}

On one hand, the extra diagonal coefficients of $(\sigma'_n)_*$ are equal to $0$ because, as $\Sigma'_n$ is an even function on $Y$, we have, for $i = 1,2$, \begin{eqnarray*}\; \left\{\!\! \begin{array}{l }
y_i \longmapsto W_n^{e_i}(y) \quad \text{is an odd function}, \\
y_i \longmapsto W_n^{e_j}(y) \quad \text{is an even function for} \ i \neq j,
\end{array} \right. \;\;\; \end{eqnarray*} which implies that $y_1 \longmapsto \Sigma'_n \nabla W_n^{e_1} \cdot \nabla W_n^{e_2}$ is an odd function. Then, by symmetry of $Y$ with respect to $0$,\[(\sigma'_n)_* e_i \cdot e_j = \int_{Y} \Sigma'_n \nabla W_n^{e_i} \cdot \nabla W_n^{e_j} \ \textrm{d}y = 0.\]

On the other hand, as $\Sigma'_n$ is isotropic, for the diagonal coefficients, we use the Voigt-Reuss inequalities (see, e.g., \cite{Jikov} p.44 or \cite{ital}): for any $i = 1,2$ and $j \neq i$, \begin{equation}
\big\langle \langle (\Sigma'_n e_i \cdot e_i)^{-1} \rangle^{-1}_i \big\rangle_j \leq (\sigma'_n)_*e_i \cdot e_i
\leq \big\langle \langle \Sigma'_n e_i \cdot e_i \rangle^{-1}_j \big\rangle_i^{-1} \label{Jik}
\end{equation} where $\langle \cdot \rangle_i$ denotes the average with respect to $y_i$ at a fixed $y_j$ for $j \neq i$.

An easy computation gives, for the direction $e_1$,\[\displaystyle (1- 2 t_n) \left( \frac{\ell - 2 t_n}{\ell \alpha'_{1,n}(h)} + \frac{2 t_n}{\ell \alpha'_{2,n}(h)} \right)^{-1} + 2 t_n \left( \frac{\ell}{\ell \alpha'_{2,n}(h)}\right)^{-1} \leq (\sigma'_n)_*e_1 \cdot e_1\] and \[ (\sigma'_n)_*e_1 \cdot e_1\leq \ell\left(\frac{\ell- 2 t_n} {(1-2 t_n) \alpha'_{1,n}(h) + 2 t_n \alpha'_{2,n}(h)} + \frac{2 t_n}{\alpha'_{2,n}(h)} \right)^{-1}.\]

\noindent By \eqref{condalphasim!} and \eqref{55555}, we have the convergence\[ \lim \limits_{n \to \infty} (\sigma'_n)_*e_1 \cdot e_1 = \alpha_1 + \frac{\alpha^2_2 + \beta^2_2 h^2}{( \ell + 1) \alpha_2}.\]

\noindent A similar computation on the direction $e_2$ gives the asymptotic behavior: \begin{equation} \displaystyle \lim \limits_{n \to \infty} (\sigma'_n)_* = \lim \limits_{n \to \infty} \left(\sigma^0_n\right)_*\big(\alpha'_{1,n}(h),\alpha'_{2,n}(h)\big) = {\begin{pmatrix} \displaystyle{\alpha_1 + \frac{\alpha^2_2 + \beta^2_2 h^2}{( \ell + 1) \alpha_2}}  & 0 \\ 0 & \displaystyle{\alpha_1+ \frac{\alpha^2_2 + \beta^2_2 h^2}{\ell ( \ell + 1) \alpha_2}} \end{pmatrix}} \label{y1!}. \end{equation}

\bigskip

Moreover, the matrix $\sigma_n(h)$ clearly satisfies all the hypothesis of Theorem \ref{Mil}. By Theorem \ref{Mil} and \eqref{y1!}, we have\begin{equation*}
\lim \limits_{n \to \infty} (\sigma_n)_*(h) = \lim \limits_{n \to \infty} \left(\sigma^0_n\right)_*\big(\alpha'_{1,n}(h),\alpha'_{2,n}(h)\big) + \beta_1 h J = {\begin{pmatrix} \displaystyle{\alpha_1 + \frac{\alpha^2_2 + \beta^2_2 h^2}{( \ell + 1) \alpha_2}}  & - h \beta_1 \\ h \beta_1 & \displaystyle{\alpha_1+ \frac{\alpha^2_2 + \beta^2_2 h^2}{\ell ( \ell + 1) \alpha_2}} \end{pmatrix}}.
\end{equation*}

\noindent We finally apply Theorem \ref{thunifconv} to get that $\sigma_*(h) = \lim \limits_{n \to \infty} (\sigma_n)_*(h)$. \qed

\section{A three-dimensional fibered microstructure\label{dim3}}

\begin{figure}[H]
\centering
\leavevmode
\subfloat[The fibers lattice]{\label{figdim3gene}\begin{tikzpicture}[math3d, scale=3]
\coordinate (a) at (1,0,0);
\coordinate (b) at (1,1,0);
\coordinate (c) at (0,1,0);
\coordinate (d) at (0,0,0);
\coordinate (e) at (1,0,1);
\coordinate (f) at (1,1,1);
\coordinate (g) at (0,1,1);
\coordinate (h) at (0,0,1);
\coordinate (a1) at (1,0,1);
\coordinate (b1) at (0.66,0,1);
\coordinate (c1) at (0.66,0.33,1);
\coordinate (d1) at (1,0.33,1);
\coordinate (e1) at (1,0,0.66);
\coordinate (f1) at (0.66,0,0.66);
\coordinate (g1) at (0.66,0.33,0.66);
\coordinate (h1) at (1,0.33,0.66);
\foreach \col in {1/6-0.07,3/6-0.07,5/6-0.07}{\foreach \h in {1/6,3/6,5/6}{
\fill[fill=gray!20, opacity=0.7] (\h,\col,0) -- (\h,\col,1)-- (\h,\col+0.14,1) -- (\h,\col+0.14,0) -- cycle;}}
\foreach \col in {1/6,3/6,5/6}{\foreach \h in {1/6,3/6,5/6}{
\draw[color=white,domain=0:360,samples=100,fill=gray!20, opacity = 0.7
] plot ({\col + 0.07*cos(\x)},{\h + 0.07*sin(\x)},{0});}}
\foreach \col in {1/6,3/6,5/6}{\foreach \h in {1/6,3/6,5/6}{
\draw[color=white, domain=0:360,samples=100,fill=gray!20, opacity = 0.7] plot ({\col + 0.07*cos(\x)},{\h + 0.07*sin(\x)},{1});}}
\foreach \col in {1/6-0.07,3/6-0.07,5/6-0.07}{\foreach \h in {1/6,3/6,5/6}{
\draw (\h,\col,0) -- (\h,\col,1);}}
\foreach \col in {1/6-0.07,3/6-0.07,5/6-0.07}{\foreach \h in {1/6,3/6,5/6}{
\draw (\h,\col+0.14,0) -- (\h,\col+0.14,1);}}
\foreach \col in {1/6,3/6,5/6}{\foreach \h in {1/6,3/6,5/6}{
\draw[domain=0:360,samples=100
] plot ({\col + 0.07*cos(\x)},{\h + 0.07*sin(\x)},{0});}}
\foreach \col in {1/6,3/6,5/6}{\foreach \h in {1/6,3/6,5/6}{
\draw[domain=0:360,samples=100] plot ({\col + 0.07*cos(\x)},{\h + 0.07*sin(\x)},{1});}}

\fill[fill=gray!20, opacity=0.7] plot ({1/6 + 0.07*cos(\x)},{1/6 + 0.07*sin(\x)},{0}) -- (0,0,0) -- cycle;
\draw[dashed] (e1) -- (f1) -- (g1) -- (h1) --cycle;
\draw[dashed] (a1) -- (e1) (b1) -- (f1) --cycle;
\draw[dashed] (c1) -- (g1) (d1) -- (h1) --cycle;
\draw[dashed] (b1) -- (c1) (c1) -- (d1);
\draw (a) -- (b) -- (c) -- (d) -- cycle;
\draw (e) -- (f) -- (g) -- (h) -- cycle;
\draw (a) -- (e) ;
\draw (b) -- (f) ;
\draw (c) -- (g) ;
\draw (d) -- (h) ;

\end{tikzpicture}}

\subfloat[The period cell]{\label{figdim3cell}\begin{tikzpicture}[math3d, scale=3]
\coordinate (a) at (1,0,0);
\coordinate (b) at (1,1,0);
\coordinate (c) at (0,1,0);
\coordinate (d) at (0,0,0);
\coordinate (e) at (1,0,1);
\coordinate (f) at (1,1,1);
\coordinate (g) at (0,1,1);
\coordinate (h) at (0,0,1);
\coordinate (cc) at (0.5,0.5,0.5);
\fill [ color=gray!20, opacity = 0.7] (0.463,0.593,1) -- (0.463,0.593,0) -- (0.537,0.407,0) -- (0.537,0.407,1)  --cycle;
\draw[color=white, domain=0:360,samples=100,fill=gray!20, opacity = 0.7] plot ({1/2 + 0.1*cos(\x)},{1/2 + 0.1*sin(\x)},{1});
\draw[color=white, domain=0:360,samples=100,fill=gray!20, opacity = 0.7] plot ({1/2 + 0.1*cos(\x)},{1/2 + 0.1*sin(\x)},{0});
\draw (0.463,0.593,1) -- (0.463,0.593,0) -- cycle;
\draw (0.537,0.407,1) -- (0.537,0.407,0) --cycle;
\draw[domain=0:360,samples=100] plot ({1/2 + 0.1*cos(\x)},{1/2 + 0.1*sin(\x)},{1});
\draw[domain=0:360,samples=100] plot ({1/2 + 0.1*cos(\x)},{1/2 + 0.1*sin(\x)},{0});
\draw[<->] (1,-0.2,1) -- (1,-0.2,0)
node[midway, left]{$\varepsilon_n$};
\draw (cc) circle (0.001);
\draw[] (d) edge (a) edge (c) edge (h);
\draw (e) -- (h) -- (g) (b)--(f);
\draw (a) -- (e) -- (f) -- (g) -- (c) -- (b) -- cycle; 
\draw (e) -- (h) -- (g) (b)--(f);
\draw[dashed] (1,-0.2,1) -- (1,0,1);
\draw[dashed] (1,-0.2,0) -- (1,0,0);
\draw[<->] (0.5,0.5,1.2) -- (0.463,0.593,1.2)
node[midway, above]{$\varepsilon_n r_n$};
\draw[] (0.5,0.5,1) -- (0.463,0.593,1) ;
\draw[dashed] (0.5,0.5,1.2) -- (0.5,0.5,1) ;
\draw[dashed] (0.463,0.593,1.2) -- (0.463,0.593,1) ;
\fill [ color=gray!20, opacity = 0.7]
(0,1.3,1) -- (0,1.4,1) -- (0,1.4,0.9) -- (0,1.3,0.9) -- cycle;
\draw (0,1.3,1) -- (0,1.4,1) -- (0,1.4,0.9) -- (0,1.3,0.9) -- cycle;
\draw (0,-1.4,0.95);
\draw (0,1.4,0.95) node[right]{$\Omega_n \cap \varepsilon_n(Y+\nu)$};
\draw[->] (0,1.35,0) -- (0.33,1.35,0)
node[below]{$x_1$};
\draw[->] (0,1.35,0) -- (0,1.68,0)
node[below]{$x_2$};
\draw[->] (0,1.35,0) -- (0,1.35,0.33)
node[above]{$x_3$};
\end{tikzpicture}}
\caption{The fibered structure in dimension 3}
\label{fiber}
\end{figure}
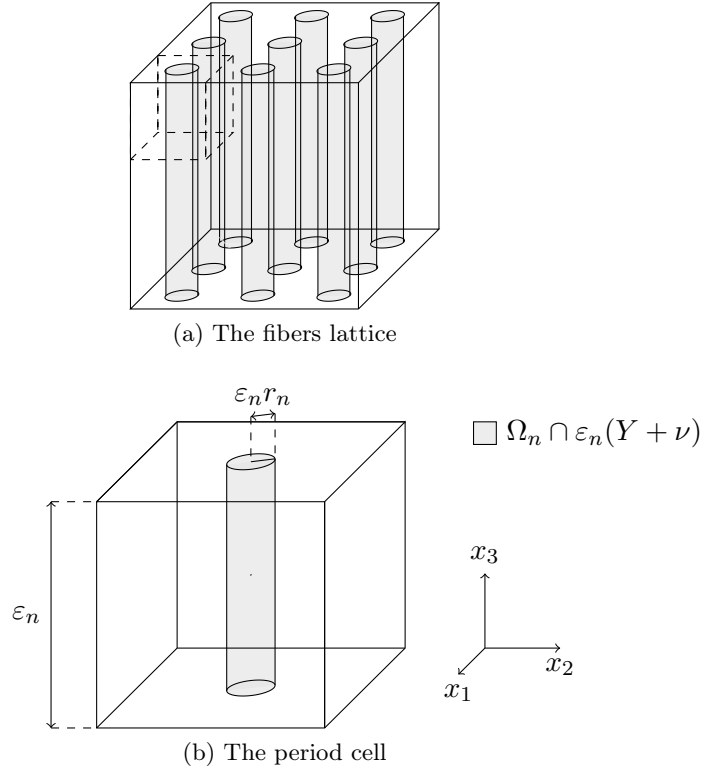

In this section we study a particular two-phase composite in dimension three. One of the phases is composed by a periodic set of high conductivity fibers embedded in an isotropic medium (figure \ref{figdim3gene}). The conductivity $\sigma_n(h)$ is not symmetric due to the perturbation of a magnetic field.

First, describe the geometry of the microstructure. Let $Y := \left(-\frac{1}{2} , \frac{1}{2} \right)^3$ be the unit cube centered at the origin of $\R^3$. For $r_n\in\left(0,{\frac{1}{2}}\right)$, consider the closed cylinder $\omega_n$ parallel to the $x_3$-axis, of radius $r_n$ and centered in $Y$: \begin{equation} \omega_n := \left\{y\in Y \ | \ y_1^2 + y_2^2 \leq r_n^2 \right\}.\end{equation} Let $\Omega=\widetilde{\Omega}\times(0,1)$ be an open cylinder of $\R^3$, where $\widetilde{\Omega}$ is a bounded domain of $\R^2$ with a Lipschitz boundary. For $\varepsilon_n\in\left(0,1\right)$, consider the closed subset $\Omega_n$ of $\Omega$ defined by the intersection with $\Omega$ of the $\varepsilon_n Y$-periodic network in $\R^3$ composed by the closed cylinders parallel to the $x_3$-axis, centered on the points $\varepsilon_n k$, $k\in\Z^2$, in the $x_1$-$x_2$ plane, and of radius $\varepsilon_n r_n$, namely: \begin{equation} \Omega_n := \Omega \cap \bigcup_{\nu \in \mathbb Z^3} \varepsilon_n(\omega_n + \nu). \end{equation} The period cell of the microstructure is represented in figure \ref{figdim3cell}.

\noindent We then define the two-phase conductivity by

\begin{eqnarray}\mathcal{\sigma}_n(h) = \; \left\{\!\! \begin{array}{c l}
\alpha_1 I_3 + \beta_1 \mathscr{E}(h) & \ \ \mathrm{in} \ \Omega\setminus \Omega_n\\
\alpha_{2,n} I_3 + \beta_{2,n} \mathscr E(h) & \ \ \mathrm{in} \ \Omega_n,
\end{array} \right. \;\;\; \label{sigmadim3}\end{eqnarray}

\noindent where $\alpha_1>0$, $\beta_1\in\R$, $\alpha_{2,n}>0$ and $\beta_{2,n}$ are real sequences, and
\[\mathscr{E}(h) := \begin{pmatrix} 0 & -h_3 & h_2 \\ h_3 & 0 & -h_1 \\ -h_2 & h_1 & 0 \end{pmatrix}, \quad \text{for} \ h = \begin{pmatrix} h_1 \\ h_2 \\h_3 \end{pmatrix} \in \R^3.\]

\par
Our aim is to study the homogenization problem

\begin{eqnarray}\mathcal{P}_{\Omega,n} \; \left\{\!\! \begin{array}{r c l l}
- \mathrm{div} (\sigma_n(h) \nabla u_n ) & = & f & \ \ \mathrm{in} \ \Omega\\
u_n & = & 0 & \ \ \mathrm{on} \ \partial \Omega.
\end{array} \right. \;\;\; \label{pbdim3}\end{eqnarray}

\begin{Theo}

Let $\alpha_1>0$, $\beta_1\in\R$, and let $\varepsilon_n, r_n, \alpha_{2,n}, \beta_{2,n}$, $n\in \mathbb N$, be real sequences such that $\varepsilon_n, r_n>0$ converge to $0$, $\alpha_{2,n}>0$, and
\begin{equation}\lim \limits_{n \to \infty} \varepsilon_n^2 |\ln r_n | = 0,\quad\lim \limits_{n \to \infty} |\omega_n|\,\alpha_{2,n} = \alpha_2 > 0,\quad
\lim \limits_{n \to \infty} |\omega_n|\,\beta_{2,n}= \beta_2\in\R.\label{conddim}\end{equation}
Consider, for $h \in \R^3$, the conductivity $\sigma_n(h)$ defined by (\ref{sigmadim3}).

\noindent Then, there exists a subsequence of n, still denoted by n, such that, for any $f \in H^{-1}(\Omega)$ and any $h \in \R^3$, the solution $u_n$ of $\mathcal P_{\Omega,n}$ converges weakly in $H_0^1(\Omega)$ to the solution $u$ of \begin{eqnarray}\mathcal{P}_{\Omega,*} \; \left\{\!\! \begin{array}{r c l l}
- \mathrm{div} (\sigma_*(h) \nabla u_n ) & = & f & \ \ \mathrm{in} \ \Omega\\
u & = & 0 & \ \ \mathrm{on} \ \partial \Omega,
\end{array} \right. \;\;\; \end{eqnarray} where
$\sigma_*(h)$ is given by \begin{equation}\sigma_*(h) = \alpha_1 I_3 + \displaystyle \left(\frac{\alpha_2^3 + \alpha_2 \beta_2^2 |h|^2}{\alpha_2^2 + \beta_2^2 h_3^2} \right) e_3 \otimes e_3 + \beta_1 \mathscr{E}(h) .
\label{sdim3}
\end{equation}
\label{thdim3}
\end{Theo}
\begin{Rem}
Theorem \ref{thdim3} can be actually extended to fibers with a more general cross-section. More precisely, we can replace the disk $r_n D$ of radius $r_n$ by the homothetic $r_n Q$ of any connected open set $Q$ included in the unit disk $D$, such that the present fiber $\omega_n$ is replaced by the new fiber $r_n Q\times\left(-\frac{1}{2},\frac{1}{2}\right)$ in the period cell of the microstructure.
\par
On the one hand, this change allows us to use the same test function $v_n$ \eqref{vndef} defined in the proof of Theorem \ref{thdim3}, since $v_n$ remains equal to $1$ in the new fibers due to the inclusion $Q\subset D$. On the other hand, Lemma \ref{esti} allows us to replace the disk $D$ by the open set $Q \subset D$.\label{extthdim3}
\end{Rem}
\begin{Rem}
We can also extend the result of Theorem \ref{thdim3} to an isotropic fibered microstructure composed by three similar periodic fibers lattices arranged in the three orthogonal directions $e_1,e_2,e_3$, namely \begin{equation*}
\omega_n := \bigcup_{j = 1}^3 \Big\{y\in Y \ | \ \Sum_{i \neq j} y_i^2 \leq r_n^2 \Big\} \quad \text{and} \quad \Omega_n := \Omega \cap \bigcup_{\nu \in \mathbb Z^3} \varepsilon_n(\omega_n + \nu),
\end{equation*} as represented in figure \ref{fiber2}.
Then, we derive the following homogenization conductivity:$$\sigma_*(h) = \alpha_1 I_3  + \displaystyle \Sum_{i=1}^3 \displaystyle \left(\frac{\alpha_2^3 + \alpha_2 \beta_2^2 |h|^2}{\alpha_2^2 + \beta_2^2 h_i^2} \right) e_i \otimes e_i + \beta_1 \mathscr{E}(h).$$
\begin{figure}[H]
\centering
\begin{tikzpicture}[math3d, scale=3]
\coordinate (a) at (1,0,0);
\coordinate (b) at (1,1,0);
\coordinate (c) at (0,1,0);
\coordinate (d) at (0,0,0);
\coordinate (e) at (1,0,1);
\coordinate (f) at (1,1,1);
\coordinate (g) at (0,1,1);
\coordinate (h) at (0,0,1);
\coordinate (cc) at (0.5,0.5,0.5);
\fill [ color=gray!20, opacity = 0.7] (0.463,0.593,1) -- (0.463,0.593,0) -- (0.537,0.407,0) -- (0.537,0.407,1)  --cycle;
\fill [ color=gray!20, opacity = 0.7] (0.463,1,0.593) -- (0.463,0,0.593) -- (0.537,0,0.407) -- (0.537,1,0.407)  --cycle;
\draw[color=white, domain=0:360,samples=100,fill=gray!20, opacity = 0.7] plot ({1/2 + 0.1*cos(\x)},{1/2 + 0.1*sin(\x)},{1});
\draw[color=white, domain=0:360,samples=100,fill=gray!20, opacity = 0.7] plot ({1/2 + 0.1*cos(\x)},{1/2 + 0.1*sin(\x)},{0});
\draw[color=white, domain=0:360,samples=100,fill=gray!20, opacity = 0.7] plot ({1/2 + 0.1*cos(\x)},{1},{1/2 + 0.1*sin(\x)});
\draw[color=white, domain=0:360,samples=100,fill=gray!20, opacity = 0.7] plot ({1/2 + 0.1*cos(\x)},{0},{1/2 + 0.1*sin(\x)});
\draw[color=white, domain=0:360,samples=100,fill=gray!20, opacity = 0.7] plot ({1},{1/2 + 0.1*cos(\x)},{1/2 + 0.1*sin(\x)});
\draw[color=white, domain=0:360,samples=100,fill=gray!20, opacity = 0.7] plot ({0},{1/2 + 0.1*cos(\x)},{1/2 + 0.1*sin(\x)});
\draw (0.463,0.593,1) -- (0.463,0.593,0) -- cycle;
\draw (0.537,0.407,1) -- (0.537,0.407,0) --cycle;
\draw[domain=0:360,samples=100] plot ({1/2 + 0.1*cos(\x)},{1/2 + 0.1*sin(\x)},{1});
\draw[domain=0:360,samples=100] plot ({1/2 + 0.1*cos(\x)},{1/2 + 0.1*sin(\x)},{0});
\draw[domain=0:360,samples=100] plot ({1/2 + 0.1*cos(\x)},{0},{1/2 + 0.1*sin(\x)});
\draw[domain=0:360,samples=100] plot ({1/2 + 0.1*cos(\x)},{1},{1/2 + 0.1*sin(\x)});
\draw[domain=0:360,samples=100] plot ({0},{1/2 + 0.1*cos(\x)},{1/2 + 0.1*sin(\x)});
\draw[domain=0:360,samples=100] plot ({1},{1/2 + 0.1*cos(\x)},{1/2 + 0.1*sin(\x)});
\draw (0.537,1,0.407) -- (0.537,0,0.407) --cycle;
\draw (0.463,1,0.593) -- (0.463,0,0.593) -- cycle;
\draw (1,0.442,0.582) -- (0,0.442,0.582) -- cycle;
\draw (1,0.550,0.410) -- (0,0.550,0.410) --cycle;
\draw[<->] (1,-0.2,1) -- (1,-0.2,0)
node[midway, left]{$\varepsilon_n$};
\draw (a) -- (e) -- (f) -- (g) -- (c) -- (b) -- cycle; 
\draw (e) -- (h) -- (g) (b)--(f); 
\draw[dashed] (1,-0.2,1) -- (1,0,1);
\draw[dashed] (1,-0.2,0) -- (1,0,0);
\draw[<->] (0.5,0.5,1.2) -- (0.463,0.593,1.2)
node[midway, above]{$\varepsilon_n r_n$};
\draw[] (0.5,0.5,1) -- (0.463,0.593,1) ;
\draw[dashed] (0.5,0.5,1.2) -- (0.5,0.5,1) ;
\draw[dashed] (0.463,0.593,1.2) -- (0.463,0.593,1) ;
\fill [ color=gray!20, opacity = 0.7]
(0,1.3,1) -- (0,1.4,1) -- (0,1.4,0.9) -- (0,1.3,0.9) -- cycle;
\draw (0,1.3,1) -- (0,1.4,1) -- (0,1.4,0.9) -- (0,1.3,0.9) -- cycle;
\draw (0,-1.4,0.95);
\draw (0,1.4,0.95) node[right]{$\Omega_n \cap \varepsilon_n(Y+\nu)$};
\draw[->] (0,1.35,0) -- (0.33,1.35,0)
node[below]{$x_1$};
\draw[->] (0,1.35,0) -- (0,1.68,0)
node[below]{$x_2$};
\draw[->] (0,1.35,0) -- (0,1.35,0.33)
node[above]{$x_3$};
\draw[] (d) edge (a) edge (c) edge (h);
\draw (e) -- (h) -- (g) (b)--(f);
\draw (a) -- (e) -- (f) -- (g) -- (c) -- (b) -- cycle; 
\draw (e) -- (h) -- (g) (b)--(f); 
\end{tikzpicture}
\caption{The period cell of the isotropic fibered structure in dimension 3}
\label{fiber2}
\end{figure}
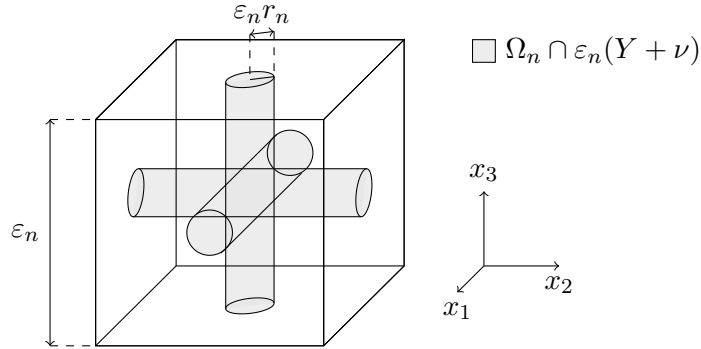
\end{Rem}

\begin{Rem}
We can check that when the volume fraction $\theta_n=\theta$ and the highly conducting phase of the conductivity $\alpha_{2,n}= \alpha_\theta$ and $\beta_{2,n}=\beta_\theta$ are independent of $n$, the explicit formula of \cite{Grab2} denoted by $\sigma_*(\theta,h)$, for the classical (since the period cell is now independent of $n$) periodically homogenized conductivity (see \eqref{defsigman*}) has a limit as $\theta \to 0$ when $\theta \alpha_\theta$ and $\theta \beta_\theta$ converge. Indeed, we may replace in the computations of \cite{Grab2} the optimal Vigdergauz shape by the circular cross-section in the previous asymptotic regime. Therefore, Theorem \ref{thdim3} validates the double process characterized by the homogenization at a fixed volume fraction $\theta$ combined with the limit as $\theta\to 0$, by one homogenization process in which both the period and the volume fraction $\theta_n = \pi r_n^2$ of the high conductivity phase tend to $0$ as $n\to\infty$.
\end{Rem}

\begin{Rem}
The hypothesis on the convergence of $\varepsilon_n^2 | \ln r_n |$ \eqref{conddim} allows us to avoid nonlocal effects in dimension three (see \cite{FeKh,BelBou}). These effects do not appear in dimension two as shown in \cite{BrARMA}. Therefore, we can make a comparison between dimension two and dimension three based on the strong field perturbation in the absence of nonlocal effects.
\label{nonloc}
\end{Rem}

\begin{Rem}\label{Rem.he3}
If $h=h_3 e_3$, the homogenized conductivity becomes
\begin{equation*}
\sigma_*(h) = \alpha_1 I_3 + \alpha_2\,e_3 \otimes e_3 + \beta_1 h_3 \left(\begin{smallmatrix} 0 & -1 & 0 \\ 1 & 0 &0 \\ 0& 0 &0\end{smallmatrix}\right)
\end{equation*}
which reduces to the simplified two-dimensional case when the symmetric part of the conductivity is independent of $h_3$ (\textit{i.e.} $\sigma_*^0$ in \eqref{defsigma*} does not depend on its second argument).
\end{Rem}

\noindent \textbf{Proof of Theorem \ref{thdim3}} The proof will be divided into four parts. We first prove the weak-$*$ convergence in $\mathcal M (\Omega)$ of $\sigma_n(h) \nabla u_n$ in $\Omega_n$. Then we establish a linear system satisfied by the limits defined by  \[\displaystyle \frac{\mathds{1}_{\Omega_n}}{|\omega_n|} \frac{\partial u_n}{\partial x_i} \harpoon 0 \xi_i \quad \text{weakly-$*$ in} \ \mathcal{M}(\Omega).\] Moreover, we deduce from Lemma \ref{esti} that \[\displaystyle \frac{\mathds{1}_{\Omega_n}}{|\omega_n|} \frac{\partial u_n}{\partial x_3} \harpoon 0 \frac{\partial u}{\partial x_3} \quad \text{weakly-$*$ in} \ \mathcal{M}(\Omega).\] We finally calculate the homogenized matrix.

\bigskip

We first remark that, classically, the sequence of solutions $u_n$ of $\mathcal{P}_{\Omega,n}$ (see (\ref{pbdim3})) is bounded in $H_0^1(\Omega)$ because, since $\alpha_{2,n}$ diverges to $\infty$ :\[ ||\nabla u_n||_{L^2(\Omega)^3}^2 \leq C \int_\Omega ( \alpha_1 \mathds{1}_{\Omega \setminus \Omega_n}I_3 + \alpha_{2,n} \mathds{1}_{ \Omega_n}I_3) \nabla u_n \cdot \nabla u_n \dx = \int_\Omega \sigma_n(h) \nabla u_n \cdot \nabla u_n \dx.\]

\noindent By the Poincaré inequality, the previous inequality and \eqref{pbdim3} lead us to \[ ||u_n||_{H^1_0(\Omega)}^2 \leq C ||\nabla u_n||_{L^2(\Omega)^3}^2 \leq C \big| \langle f , u_n \rangle_{H^{-1}(\Omega),H_0^1(\Omega)}\big| \leq C ||f||_{H^{-1}(\Omega)} ||u_n||_{H_0^{1}(\Omega)}\] and then to\[ ||u_n||_{H_0^{1}(\Omega)} \leq C ||f||_{H^{-1}(\Omega)}.\]

\noindent Thus, up to a subsequence still denoted by $n$, $u_n$ converges weakly to some function $u$ in  $H_0^1(\Omega)$.

\bigskip

\noindent \textit{First step}: Weak-$*$ convergence in $\mathcal M(\Omega)$ of the conductivity in the fibers $\mathds{1}_{\Omega_n} \ \big(\alpha_{2,n} I_3 + \beta_{2,n} \mathscr{E}(h)\big) \nabla u_n$.
\smallskip
We proceed as in \cite{BrTc} with a suitable oscillating test function. For $R \in (0,1/2)$, define the $Y$-periodic (independent of $y_3$) function $V_n$ by
\[V_n(y_1,y_2,y_3) = \; \left\{\!\! \begin{array}{c l}
1 & \ \ \mathrm{if} \ \sqrt{y_1^2+y_2^2} \leq r_n\\
\displaystyle \frac{\ln R - \ln \sqrt{y_1^2+y_2^2}}{\ln R - \ln r_n} & \ \ \mathrm{if} \ r_n \leq \sqrt{y_1^2+y_2^2} \leq R \\
0 & \ \ \mathrm{if} \ \sqrt{y_1^2+y_2^2} \geq R,
\end{array} \right. \quad\mbox{for }y\in Y,
\]
\noindent and the rescaled function \begin{equation}v_n(x) = V_n\left(\frac{x}{\varepsilon_n}\right),\quad\mbox{for }x\in\R^3.\label{vndef}\end{equation}
In particular, by using the cylindrical coordinates and the fact that $r_n$ converges to $0$, this function satisfies the inequalities\begin{align*}||v_n||^2_{L^2(\Omega)} \leq C ||V_n||^2_{L^2(Y)} &= C \left|\ln \frac{R}{r_n}\right|^{-2} \left(\pi r_n^2 + \int_0^{2 \pi} \int_{r_n}^R r \ln^2 \frac{R}{r} \ \textrm{d}r \textrm{d}\theta \right)\\
&= C \left|\ln \frac{R}{r_n}\right|^{-2} \left( \pi \frac{R^2 - r_n^2}{2} - \pi r_n^2 \ln^2 \frac{R}{r_n} - \pi \ln \frac{R}{r_n} \right) \leq C \left|\ln \frac{R}{r_n}\right|^{-2},\end{align*}\begin{equation*}\hspace{-0cm} ||\nabla v_n||^2_{L^2(\Omega)^3} \leq \displaystyle \frac{C}{\varepsilon_n^2} ||\nabla V_n||^2_{L^2(Y)^3} = \frac{C}{\varepsilon_n^2} \left|\ln \frac{R}{r_n}\right|^{-2}\int_0^{2 \pi} \int_{r_n}^R \frac{1}{r} \ \textrm{d}r \textrm{d}\theta \leq \frac{C}{\varepsilon_n^2} \left|\ln \frac{R}{r_n}\right|^{-1}\end{equation*}and, consequently\begin{equation}||v_n||_{L^2(\Omega)} + \varepsilon_n ||\nabla v_n||_{L^2(\Omega)^3} \leq C \sqrt{\left|\ln \frac{R}{r_n}\right|^{-1}} \underset{n\to \infty}{\longrightarrow} 0 . \label{Dpropvn}\end{equation}

Let $\lambda$ be a vector in $\R^3$ perpendicular to the $x_3$-axis. Define the $Y$-periodic function $\widetilde{X}_n$ by $\nabla \widetilde{X}_n = \lambda$ in $\omega_n$, such that $\widetilde{X}_n \in \mathscr{D}(Y)$ and is $Y$-periodic, and the rescaled function $X_n$ by \begin{equation} X_n(x) = \varepsilon_n \ \widetilde{X}_n\left(\frac{x}{\varepsilon_n}\right).\label{DXn}\end{equation}In particular, $X_n$ satisfies \begin{equation} ||X_n||_\infty = \varepsilon_n \ ||\widetilde{X}_n||_\infty \leq C \varepsilon_n \quad , \quad ||\nabla X_n||_\infty = ||\nabla \widetilde{X}_n||_\infty \leq C \quad \text{and} \quad \nabla X_n = \lambda \quad \text{in} \ \Omega_n . \label{DpropXn}\end{equation} We have, by \eqref{DpropXn} and \eqref{Dpropvn}, \begin{align*}||v_n X_n||_{H^1(\Omega)} &\leq ||X_n||_\infty ||v_n||_{L^2(\Omega)} + ||X_n||_\infty ||\nabla v_n||_{L^2(\Omega)^3} +||\nabla X_n||_\infty ||v_n||_{L^2(\Omega)} \\ &\leq C \big(||v_n||_{L^2(\Omega)} + \varepsilon_n ||\nabla v_n||_{L^2(\Omega)^3} \big) \underset{n\to \infty}{\longrightarrow} 0, \end{align*}which gives\begin{equation}
\forall\,\varphi \in \mathscr{D}(\Omega), \quad \varphi \ v_n X_n \underset{n\to \infty}{\longrightarrow} 0 \quad \text{strongly in} \ H_0^1(\Omega). \label{gggg}
\end{equation}

Let $\varphi \in \mathscr{D}(\Omega)$. By the strong convergence \eqref{gggg}, we have \begin{equation}\int_\Omega \sigma_n(h) \nabla u_n \cdot \nabla \big( \varphi \ v_n X_n \big)\dx = \langle f , \varphi \ v_n X_n \rangle_{H^{-1}(\Omega),H_0^1(\Omega)} \underset{n\to \infty}{\longrightarrow} 0 .\label{Dmain}\end{equation}

Let us decompose this integral which converges to 0, into the integral on the fibers set $\Omega_n$ and the integral on its complementary:
\begin{subequations}
\begin{align}
\label{eq1}\int_\Omega \sigma_n(h) \nabla u_n \cdot \nabla \big( \varphi \ v_n X_n \big) \dx & = \int_{\Omega \setminus \Omega_n} (\alpha_1 I_3+\beta_1\mathscr{E}(h)) \nabla u_n \cdot \nabla \big( \varphi \ v_n X_n \big) \dx \\
\label{eq2}&+ \int_{\Omega_n} (\alpha_{2,n} I_3 + \beta_{2,n} \mathscr{E}(h)) \nabla u_n \cdot \nabla \big( \varphi \ v_n X_n \big) \dx.
\end{align}
\end{subequations}

The expression \eqref{eq1} converges to $0$ since, by the Cauchy-Schwarz inequality, the boundedness of $u_n$ in $H_0^1(\Omega)$ and \eqref{gggg}, we have \begin{equation}
\left|\int_{\Omega \setminus \Omega_n} \hspace{-0.65cm} (\alpha_1 I_3+\beta_1 \mathscr{E}(h)) \nabla u_n \cdot \nabla \big( \varphi \ v_n X_n \big) \dx\right| \leq\left|\alpha_1 I_3+\beta_1 \mathscr{E}(h)\right| ||\nabla u_n||_{L^2(\Omega)^3} \ ||\varphi \ v_n X_n||_{H_0^1(\Omega)} \underset{n\to \infty}{\longrightarrow} 0. \label{convfibb}
\end{equation}Consequently, as $v_n = 1$ and $\nabla X_n = \lambda$ on $\Omega_n$, by \eqref{Dmain}, \eqref{eq1}, \eqref{eq2} and \eqref{convfibb}, we have\begin{equation}\int_{\Omega_n} \sigma_n(h) \nabla u_n \cdot \lambda \ \varphi \dx + \int_{\Omega_n} \sigma_n(h) \nabla u_n \cdot \nabla \varphi \ X_n \dx \underset{n\to \infty}{\longrightarrow} 0. \label{Dfiber}\end{equation}

To prove the convergence to $0$ of the right term, we now show that $ \mathds{1}_{\Omega_n} \big(\alpha_{2,n} I_3 + \beta_{2,n} \mathscr{E}(h)\big) \nabla u_n$ is bounded in $L^1(\Omega)^3$. We have, by the Cauchy-Schwarz inequality, \eqref{conddim} and the classical equivalent $|\Omega_n| \underset{n \to \infty}{\sim} |\Omega| \ |\omega_n|$,\begin{align*} \displaystyle \left(\int_{\Omega_n} \big|\big(\alpha_{2,n} I_3 + \beta_{2,n} \mathscr{E}(h)\big) \nabla u_n \big|\dx \right)^2 & \leq \big|I_3 + \alpha_{2,n}^{-1} \beta_{2,n} \mathscr{E}(h)\big|^2\,|\Omega_n|\,\alpha_{2,n}\int_{\Omega_n} \alpha_{2,n} |\nabla u_n|^2 \dx \\ 
& \leq C\int_{\Omega}\sigma_n(h) \nabla u_n \cdot \nabla u_n \dx \\
&\leq C\,||f||_{H^{-1}(\Omega)}\,||u_n||_{H_0^1(\Omega)}. \end{align*}This combined with the boundedness of $u_n$ in $H_0^1(\Omega)$ implies that $\mathds{1}_{\Omega_n} \big(\alpha_{2,n} I_3 + \beta_{2,n} \mathscr{E}(h)\big) \nabla u_n$ is bounded in $L^1(\Omega)^3$. This bound and the uniform convergence to $0$ of $X_n$ (see \eqref{DpropXn}) imply the convergence to $0$ of the right term of \eqref{Dfiber}, hence \[ \displaystyle \int_{\Omega_n} \big(\alpha_{2,n} I_3 + \beta_{2,n} \mathscr{E}(h)\big) \nabla u_n \cdot \lambda \ \varphi \dx \underset{n\to \infty}{\longrightarrow} 0.\] We rewrite this condition as \begin{equation} \forall\,\lambda \perp e_3, \ \ \mathds{1}_{\Omega_n} \ \big(\alpha_{2,n} I_3 + \beta_{2,n} \mathscr{E}(h)\big) \nabla u_n \cdot \lambda \ \harpoon 00 \quad \textrm{weakly-$*$ in} \ \mathcal{M}(\Omega). \label{bgconv}\end{equation}

\bigskip

\noindent \textit{Second step}: Linear relations between weak-$*$ limits of $\displaystyle \frac{\mathds{1}_{\Omega_n}}{|\omega_n|} \frac{\partial u_n}{\partial x_i}$.

\smallskip

\noindent Thanks to the Cauchy-Schwarz inequality, we have

\[\displaystyle \left|\left| \frac{\mathds{1}_{\Omega_n}}{|\omega_n|} \frac{\partial u_n}{\partial x_i}\right|\right|_{L^1(\Omega)} \leq \frac{1}{|\omega_n|}\int_{\Omega_n} |\nabla u_n| \dx\leq \frac{1}{\sqrt{\alpha_{2,n} |\omega_n|}} \sqrt{\frac{|\Omega_n|}{|\omega_n|}}\sqrt{\int_{\Omega_n}\alpha_{2,n} |\nabla u_n|^2 \dx}\]which leads us, by \eqref{conddim} and the asymptotic behavior $|\Omega_n| \underset{n \to \infty}{\sim} |\Omega| \ |\omega_n|$, to
$$ 
		\begin{array}{rl}
			\displaystyle \left|\left| \frac{\mathds{1}_{\Omega_n}}{|\omega_n|} \frac{\partial u_n}{\partial x_i}\right|\right|_{L^1(\Omega)} & \leq \displaystyle \frac{C}{\sqrt{\alpha_{2,n} |\omega_n|}} \int_{\Omega} \sigma_n(h) \nabla u_n \cdot \nabla u_n \dx \leq C \left|\langle f , u_n \rangle_{H^{-1}(\Omega),H_0^1(\Omega)} \right|
		\end{array}$$ which is bounded by the boundedness of $u_n$ in $H_0^1(\Omega)$. This allows us to define, up to a subsequence, the following limits \begin{equation}\displaystyle \frac{\mathds{1}_{\Omega_n}}{|\omega_n|} \frac{\partial u_n}{\partial x_i} \rightharpoonup \xi_i \quad \textrm{weakly-$*$ in} \ \mathcal{M}(\Omega), \quad \text{for} \ i=1,2,3.\label{Ddefxi}\end{equation}

\noindent Then, by \eqref{bgconv} we have  \[\displaystyle \big(\alpha_{2,n} I_3 + \beta_{2,n} \mathscr{E}(h)\big)\mathds{1}_{\Omega_n} \nabla u_n \cdot \lambda = \big(\alpha_{2,n} |\omega_n| I_3 + \beta_{2,n} |\omega_n| \mathscr{E}(h)\big) \frac{\mathds{1}_{\Omega_n}}{|\omega_n|} \nabla u_n \cdot \lambda \rightharpoonup 0 \ \ \textrm{weakly-$*$ in} \ \mathcal{M}(\Omega).\]
Therefore, putting $\lambda = e_1,e_2$ in this limit and using condition (\ref{conddim}), we obtain the linear system\begin{equation*}
  \left\{
      \begin{aligned}
        \alpha_2 \xi_1 + \beta_2 h_2 \xi_3 - \beta_2 h_3 \xi_2 &= 0 \\
\alpha_2 \xi_2 + \beta_2 h_3 \xi_1 - \beta_2 h_1 \xi_3 &= 0
      \end{aligned}
    \right. \quad \text{in} \ \mathcal M(\Omega),
\end{equation*}
 which is equivalent to\begin{equation}
  \left\{
      \begin{aligned}
        \displaystyle \xi_1 &= \frac{\beta_2^2 h_1 h_3 - \alpha_2 \beta_2 h_2}{\alpha_2^2 + \beta_2^2 h_3^2} \xi_3 \\
\xi_2 &= \frac{\beta_2^2 h_2 h_3 + \alpha_2 \beta_2 h_1}{\alpha_2^2 + \beta_2^2 h_3^2} \xi_3
      \end{aligned}
    \right. \quad \text{in} \ \mathcal M(\Omega).
    \label{lin}
\end{equation}

\bigskip

\noindent \textit{Third step}: Proof of $\xi_3 = \displaystyle \frac{\partial u}{\partial x_3}$.

\smallskip

\noindent We need the following result which is an extension of the estimate (3.13) of \cite{BriMokMur}. The statement of this lemma is more general than necessary for our purpose but is linked to Remark \ref{extthdim3}.

\begin{Lem} \label{esti} Let $Q$ be a non-empty connected open subset of the unit disk $D$. Then, there exists a constant $C >0$ such that any function $U \in H^1(Y)$ satisfies the estimate\begin{equation} \displaystyle \left|\frac{1}{|r_n Q|} \int_{r_n Q \times \left( -\frac{1}{2} , \frac{1}{2}\right)} U \ \mathrm{d}y - \int_{Y} U \ \mathrm{d}y \right| \leq C \sqrt{|\ln r_n|} \ ||\nabla U||_{L^2(Y)^3}.\label{ineqgrad1}\end{equation} \label{ineqgrad}
\end{Lem}

\noindent \textbf{Proof of Lemma \ref{ineqgrad}}. Let $U \in H^1(Y)$. To prove Lemma \ref{esti}, we compare the average value of $U$ on $r_n Q$ and $r_n D$. Denoting $\widetilde{y} = (y_1,y_2)$, we have, for any $y_3 \in \left( -\frac{1}{2} , \frac{1}{2}\right)$, \begin{align*} \left|\fint_{r_n Q} U(\widetilde{y},y_3) \ \mathrm{d}\widetilde{y} - \fint_{r_n D} U(\widetilde{y},y_3) \ \mathrm{d}\widetilde{y} \right| &= \left|\fint_{Q} U(r_n \widetilde{y},y_3) \ \mathrm{d}\widetilde{y} - \fint_{D} U(r_n \widetilde{y},y_3) \ \mathrm{d}\widetilde{y} \right| \\
&\leq \fint_{Q} \left| U(r_n \widetilde{y},y_3)  - \fint_{D} U(r_n \widetilde{y},y_3) \ \mathrm{d}\widetilde{y} \right| \ \mathrm{d}\widetilde{y} , \end{align*} and, since $Q \subset D$, \begin{align*} \left|\fint_{r_n Q} U(\widetilde{y},y_3) \ \mathrm{d}\widetilde{y} - \fint_{r_n D} U(\widetilde{y},y_3) \ \mathrm{d}\widetilde{y} \right| 
&\leq \frac{|D|}{|Q|} \fint_{D} \left| U(r_n \widetilde{y},y_3)  - \fint_{D} U(r_n \widetilde{y},y_3) \ \mathrm{d}\widetilde{y} \right| \ \mathrm{d}\widetilde{y} \\
& \leq C\fint_D r_n \left(\left|\frac{\partial U}{\partial x_1}\right| + \left|\frac{\partial U}{\partial x_2}\right|\right)(r_n \widetilde{y},y_3) \ \mathrm{d}\widetilde{y} \\
&= \frac{C}{\pi r_n} \int_{r_n D} \left(\left|\frac{\partial U}{\partial x_1}\right| + \left|\frac{\partial U}{\partial x_2}\right|\right)( \widetilde{y},y_3) \ \mathrm{d}\widetilde{y} , \end{align*} the last inequality being a consequence of the Poincaré-Wirtinger inequality. Hence, integrating the previous inequality with respect to $y_3 \in \left( -\frac{1}{2} , \frac{1}{2}\right)$ and applying the Cauchy-Schwarz inequality, we obtain that \begin{align*}\left|\fint_{r_n Q \times \left( -\frac{1}{2} , \frac{1}{2}\right)} U(y) \ \mathrm{d}y - \fint_{r_n D \times \left( -\frac{1}{2} , \frac{1}{2}\right)} U(y) \ \mathrm{d}y \right| &\leq \frac{C}{\pi r_n} \int_{r_n D \times \left( -\frac{1}{2} , \frac{1}{2}\right)} \left|\nabla U\right|(y) \ \mathrm{d}y \\
& \leq C \sqrt{\int_{r_n D \times \left( -\frac{1}{2} , \frac{1}{2}\right)} |\nabla U|^2(y) \ \mathrm{d}y} \\ 
& \leq C  ||\nabla U||_{L^2(Y)^3}.\end{align*} This combined with the estimate (3.13) of \cite{BriMokMur}, i.e. \eqref{ineqgrad1} for $Q = D$, and the fact that $\sqrt{ | \ln r_n|}$ diverges to $\infty$ give the thesis. \qed

\bigskip

Let $\varphi \in \mathscr D(\Omega)$. A rescaling of \eqref{ineqgrad1} with $Q= D$ implies the inequality \begin{align*} \displaystyle \left|\frac{1}{|\omega_n|} \int_{\Omega_n} u_n \varphi \dx - \int_{\Omega} u_n \varphi \dx \right| \leq C \varepsilon_n \sqrt{|\ln r_n|} \displaystyle \left|\left|\nabla (u_n \varphi) \right|\right|_{L^2(\Omega)^3}.
\end{align*}
Combining this estimate and the first condition of \eqref{conddim} with
\[\left|\left|\nabla (u_n \varphi) \right|\right|_{L^2(\Omega)^3} \leq \left|\left|\nabla u_n \right|\right|_{L^2(\Omega)^3}||\varphi ||_\infty + \left|\left|u_n \right|\right|_{L^2(\Omega)}||\nabla \varphi ||_\infty \leq C, \] it follows that \[\displaystyle \frac{\mathds{1}_{\Omega_n}}{|\omega_n|} u_n - u_n \rightharpoonup 0 \ \ \textrm{in} \ \mathscr{D}'(\Omega).\]This convergence does not hold true when $\varepsilon_n^2| \ln r_n|$ converges to some positive constant. Under this critical regime, non-local effects appear (see Remark \ref{nonloc}). 

\par
Finally, as $\mathds{1}_{\Omega_n}$ does not depend on the $x_3$ variable, we have\[\frac{\mathds{1}_{\Omega_n}}{|\omega_n|} \frac{\partial u_n}{\partial x_3} = \frac{\partial}{\partial x_3}\frac{\mathds{1}_{\Omega_n}}{|\omega_n|}u_n = \frac{\partial}{\partial x_3}\left(\frac{\mathds{1}_{\Omega_n}}{|\omega_n|}u_n - u_n\right) + \frac{\partial u_n}{\partial x_3} \rightharpoonup \frac{\partial u}{\partial x_3} = \xi_3 \ \ \textrm{in} \ \mathscr{D}'(\Omega).\]

\bigskip

\noindent\textit{Fourth step}: Derivation of the homogenized matrix.

\smallskip

\noindent We now study the limit of $\sigma_n(h) \nabla u_n$ in order to obtain $\sigma_*(h)$. We have
\begin{equation} 
		\begin{array}{rl}
			\displaystyle \sigma_n(h) \nabla u_n \cdot e_1 &=\vspace{0.5cm}\displaystyle \mathds{1}_{\Omega \setminus \Omega_n}\left( \alpha_1 \frac{\partial u_n}{\partial x_1} - \beta_1 h_3 \frac{\partial u_n}{\partial x_2} + \beta_1 h_2 \frac{\partial u_n}{\partial x_3}\right) \\ 
			 &  \displaystyle \hspace{0.5cm} + \alpha_{2,n} |\omega_n| \frac{\mathds{1}_{\Omega_n}}{|\omega_n|} \frac{\partial u_n}{\partial x_1} - \beta_{2,n} h_3 |\omega_n| \frac{\mathds{1}_{\Omega_n}}{|\omega_n|}\frac{\partial u_n}{\partial x_2}+ \beta_{2,n} h_2 |\omega_n| \frac{\mathds{1}_{\Omega_n}}{|\omega_n|}\frac{\partial u_n}{\partial x_3}.
		\end{array}\label{5555}\end{equation} Hence, passing to the weak-$*$ limit in $\mathcal M(\Omega)$ this equality and using the linear system \eqref{lin}, $\sigma_n(h) \nabla u_n \cdot e_1$ weakly-$*$ converges in $\mathcal M(\Omega)$ to\begin{align*}
&\left( \alpha_1 \frac{\partial u}{\partial x_1} - \beta_1 h_3 \frac{\partial u}{\partial x_2} + \beta_1 h_2 \frac{\partial u}{\partial x_3} \right) + \alpha_2 \xi_1 - \beta_2 h_3 \xi_2 + \beta_2 h_2 \xi_3 \\
&=\big(\alpha_1 I_3 + \beta_1 \mathscr{E}(h)\big)\nabla u \cdot e_1 + \alpha_2 \frac{\beta_2^2 h_1 h_3 - \alpha_2 \beta_2 h_2}{\alpha_2^2 + \beta_2^2 h_3^2} \xi_3 - \beta_2 h_3 \frac{\beta_2^2 h_2 h_3 + \alpha_2 \beta_2 h_1}{\alpha_2^2 + \beta_2^2 h_3^2} \xi_3 + \beta_2 h_2 \xi_3 \\
&=\big(\alpha_1 I_3 + \beta_1 \mathscr{E}(h)\big)\nabla u \cdot e_1 \\
&\hspace{0.5cm} + \underbrace{\frac{\alpha_2(\beta_2^2 h_1 h_3 - \alpha_2 \beta_2 h_2) - \beta_2 h_3 (\beta_2^2 h_2 h_3 + \alpha_2 \beta_2 h_1) + \beta_2 h_2 (\alpha_2^2 + \beta_2^2 h_3^2)}{\alpha_2^2 + \beta_2^2 h_3^2}}_{= 0} \xi_3,
\end{align*} that is \begin{equation}
\sigma_n(h) \nabla u_n \cdot e_1 \ \harpoon 0\ \big(\alpha_1 I_3 + \beta_1 \mathscr{E}(h)\big)\nabla u \cdot e_1\quad \text{weakly-$*$ in } \mathcal M(\Omega).
\label{deriv1}\end{equation} The same calculus leads us to\begin{equation}
\sigma_n(h) \nabla u_n \cdot e_2 \ \harpoon 0\ \big(\alpha_1 I_3 + \beta_1 \mathscr{E}(h)\big)\nabla u \cdot e_2\quad \text{weakly-$*$ in } \mathcal M(\Omega).
\label{deriv2}\end{equation}We have, for the last direction $e_3$,
\begin{align*}\sigma_n(h) \nabla u_n \cdot e_3 \rightharpoonup \left( \alpha_1 \frac{\partial u}{\partial x_3} - \beta_1 h_2 \frac{\partial u}{\partial x_1} + \beta_1 h_1 \frac{\partial u}{\partial x_2} \right) + \alpha_2 \xi_3 + \beta_2 h_2 \xi_1 - \beta_2 h_1 \xi_2 \ \ \textrm{weakly-$*$ in} \ \mathcal{M}(\Omega). \end{align*} Hence, again with the linear system \eqref{lin},\begin{align*}
&\left( \alpha_1 \frac{\partial u}{\partial x_3} - \beta_1 h_2 \frac{\partial u}{\partial x_1} + \beta_1 h_1 \frac{\partial u}{\partial x_2} \right) + \alpha_2 \xi_3 - \beta_2 h_2 \xi_1 + \beta_2 h_1 \xi_2 \\
&=\big(\alpha_1 I_3 + \beta_1 \mathscr{E}(h)\big)\nabla u \cdot e_3 + \alpha_2 \xi_3 - \beta_2 h_2 \frac{\beta_2^2 h_1 h_3 - \alpha_2 \beta_2 h_2}{\alpha_2^2 + \beta_2^2 h_3^2} \xi_3 + \beta_2 h_1 \frac{\beta_2^2 h_2 h_3 + \alpha_2 \beta_2 h_1}{\alpha_2^2 + \beta_2^2 h_3^2} \xi_3. \end{align*}

\noindent Finally, by the previous equality, \eqref{deriv1} and \eqref{deriv2}, we get that$$\sigma_*(h) = \alpha_1 I_3  + \displaystyle \displaystyle \left(\frac{\alpha_2^3 + \alpha_2 \beta_2^2 |h|^2}{\alpha_2^2 + \beta_2^2 h_3^2} \right) e_3 \otimes e_3 + \beta_1 \mathscr{E}(h).$$ \qed

\bigskip

\noindent\textit{Acknowledgments:} The authors wish to thank the referees for suggestions which led to improvement.

\bibliographystyle{plain}
\bibliography{biblio}

\end{document}